\documentclass[11pt, sec]{article}
\usepackage{graphicx}
\usepackage{amsmath}
\usepackage{amssymb}
\usepackage{amsfonts}

\makeatletter
\@addtoreset{equation}{section}   
\makeatother

\makeatletter
\@addtoreset{figure}{section}
\def\thefigure{\thesection.\@arabic\c@figure}
\makeatother

\makeatletter
\@addtoreset{table}{section}
\def\thetable{\thesection.\@arabic\c@table}
\makeatother

\begin{document}

\title{\bf  On the Siegel-Weil Theorem for Loop Groups (II)}

\author{Howard Garland and Yongchang Zhu}
\date{}

\maketitle

 \newtheorem{theorem}{Theorem}[section]
  \newtheorem{prop}[theorem]{Proposition}
  \newtheorem{cor}[theorem]{Corollary}
  \newtheorem{lemma}[theorem]{Lemma}
  \newtheorem{defn}[theorem]{Definition}
  \newtheorem{ex}[theorem]{Example}

 \newcommand{\cx}{{\bf C}}
\newcommand{\la}{\langle}
\newcommand{\ra}{\rangle}
\newcommand{\res}{{\rm Res}}
\newcommand{\expp}{{\rm exp}}
\newcommand{\Sp}{{\rm Sp}}
\newcommand{\lgp}{\widehat{\rm G}( { F} ((t)) )}
\newcommand{\svee}{\scriptsize \vee}
\newcommand{\deff}{\stackrel{\rm def}{=}}
\maketitle

\maketitle

\section{Introduction} \label{section 1}

 This is the second of our two papers on the Siegel-Weil theorem for loop groups. In the first paper \cite{gz} we proved the Siegel-Weil 
 theorem for (finite dimensional) snt-modules (\cite{gz}, Theorem 8.1). In the present paper we use this result to obtain the Siegel-Weil
 theorem for loop groups, Theorem 7.5, below.

 In addition to the corresponding result for snt-modules, our proof depends on a convergence condition for certain Eisenstein series 
 on loop groups (Theorem 5.3, below). We note that this convergence criterion is used for the convergence 
 criterion for Eisenstein series associated with snt-modules (Theorem 6.6, below). The uniform convergence obtained in Theorem 6.6 
 is crucial for applying the abstract lemma in Weil \cite{W2} (see \cite{W2}, Proposition 2, page 7).

 The Siegel-Weil theorem for snt-modules does not immediately give the result for loop groups. The failure to do so is measured 
 by the terms on the right hand side of (7.8). However, in $\S$8, we show that in fact, these "error terms" vanish!

We now describe briefly our main result. 
  Let $F$ be a number field, $F^{2n}$ be the standard symplectic space over $F$, and let $(V, ( , ) )$ be a 
 finite dimensional $F$-space with an anisotropic non-degenerate symmetric bilinear form $(, )$ with  corresponding 
 orthogonal group $G$. The space $ F^{2n} \otimes V$ is naturally a symplectic space with isometry group $Sp_{2N} $ (where $2N = 2n {\rm dim } V$).
 The groups $Sp_{2n}$ and  $G$ are commuting subgroups
 in $Sp_{2N} $.  The Weil representation can be generalized to loop symplectic groups \cite{Z}. 
  Let 
\[ {\cal S} (  (t^{-1} F[ t^{-1} ] ^{2n} \otimes V )_{\bf A} )\]
 be the space of 
 Schwartz functions on the countably infinite dimensional adelic space 
\[  (t^{-1} F[ t^{-1} ] ^{2n} \otimes V )_{\bf A} . \]
 It 
 is a model for the Weil representation of the adelic loop  metaplectic group $ \widetilde {\rm Sp}_{2N} ( {\bf A} \la t \ra )$.
  The commuting pair  $(Sp_{2n} , G)$ is supposed to  be lifted to a commuting pair 
  consisting of the metaplectic loop group $  \widetilde {\rm Sp}_{2n} ( {\bf A} \la t \ra ) $ 
  and a central extension of $ G ( {\bf A} (( t)) )$ in $  \widetilde {\rm Sp}_{2N} ( {\bf A} \la t \ra )$.
  For our formulation of the Siegel-Weil formula, we only need half of the loop orthogonal group $G ( {\bf A} [[ t]] )$.
  We prove for functions $f\in  {\cal S} (  (t^{-1} F[ t^{-1} ] ^{2n} \otimes V )_{\bf A} )$ which satisfy certain properties, 
   that the theta functional
 \[  \theta ( f ) = \sum_{ r\in (t^{-1} F[ t^{-1} ] ^{2n} \otimes V }  f ( r ) \]
  converges (See Theorem 3.3).
 Also the Eisenstein series 
  \[   E ( f ) =    \sum_{ g \in  Sp_{2n } ( F[[ t]] ) \backslash Sp_{2n } ( F(( t))) }       ( g \cdot f) ( 0 )  \]
 converges under the condition ${\rm dim} \, V  > 6n +2 $  (Theorem 5.3). 
 In the above, $ Sp_{2n } ( F[[ t]] )$ plays the role Siegel parabolic subgroup for a loop symplectic group.
Our main result is that
\[  E ( f ) = \int_{   G ( F[[t]]) \backslash G ( {\bf A} [[t]] ) }   \theta (  g f ) dg ,\]
see Theorem 7.5, below.

  This paper is organized as follows: In $\S$2, we review the metaplectic loop group and Weil representation, as constructed in \cite{Z}.
 We also give some further constructions and technical results which will be needed later in the paper, including some discussion
 for non-archimedean local fields and for adeles. In $\S$4 we discuss Eisenstein series for loop metaplectic groups, and in $\S5$, prove 
 a convergence theorem, Theorem 5.3, for these Eisenstein series.  In $\S$6 we relate Eisenstein series in the loop case and the snt-module case, and finally in $\S$7, we prove the Siegel-Weil theorem (Theorem 7.5) for loop groups. However we give the proof that the terms on the right 
 hand side of (7.8) are zero in $\S$ 8.

 \section{ Metaplectic Loop Groups and the Weil Representation }

In this section, we recall the Weil representation of a loop symplectic group over a local field  
 constructed in \cite{Z} (Section 2.1). We then define the metaplectic loop group using the symbol of the Weil representation and 
 study its  Bruhat and Iwasawa decomposition (Section 2.2). In Section 2.3, we define 
   the adelic metaplectic loop group for a number field and  the adelic Weil representation.  

\

\noindent { \bf  2.1. The Weil representation and metapletic loop group over a local field.  }   
 Let $F$ be  a local field  of characteristic $0$
 and $\psi $ be  a non-trivial additive character of $F$. 
   For a standard $2N$-dimensional symplectic space $F^{2N}$ over $F$ with symplectic form $\la ,\ra$, the space 
 $ F((t))^{2N} = F^{2N} \otimes_F F(( t)) $ has an  
    $F((t))$-valued symplectic form $\la  , \ra_{F((t))}$ given by the scalar extension.  It gives a $F$-valued
  symplectic form on $ F((t))^{2N}$ by taking the residue:  for $ w , v \in F((t))^{2N}$, 
    \[    \la w , v \ra = \res \la w , v \ra_{F((t))}, \]   
  where $ \res \, a $ for $a \in F((t))$ is the coefficient of $t^{-1}$ in the expression of $a$.  
  The spaces $  X_- = t^{-1} F[t^{-1}]^{2N}$  and $X_+ =  F[[t]]^{2N}$ are maximal isotropic subspaces
 of $F((t))^{2N}$.  Since the group $\Sp_{2N} ( F(( t)))$ preserves $\la , \ra_{F((t))}$,
   it  preserves $\la , \ra$ on $F((t))^{2N}$.  We shall assume $\Sp_{2N} ( F(( t)))$ acts on $F((t))^{2N}$ from the right.
   We also need a larger group $ \Sp ( F((t))^{2N} , X_+ )$,  which is by definition  the group of all $F$-linear  symplectic isomorphisms $g$ of $F((t))^{2N}$
such that $ X_+ g $ and $X_+$ are commensurable.  The loop group $ \Sp_{2N} ( F(( t)))$ is a subgroup of  $ \Sp ( F((t))^{2N} , X_+)$ .

 Each $g\in \Sp ( F((t))^{2N} , X_+ )$ has a matrix form 
\begin{equation}\label{gform}
  g = \left[ \begin{array}{cc}  \alpha  &  \beta \\ \gamma & \delta  \end{array} \right] \end{equation}
 with respects to the decomposition $ F((t))^{2N} =  X_- \oplus X_+ $. So it acts on
  $v+ v^*$ as $ ( v + v^* ) g = ( v \alpha + v^* \gamma ) + ( v \beta + v^* \delta ) $,
  where $\alpha : X_-\to X_- $ , $\beta : X_-\to X_+$, $\gamma : X_+ \to X_- $ and $\delta : X_+ \to X_+ $. 
   We sometimes write $ \alpha_g, \beta_g , \gamma_g , \delta_g$ to indicate 
  the operators are associated with $g$.   For a symplectic isomorphism $g$ of $F((t))^{2N}$, the condition that  $ X_+ g $ and $X_+$ are commensurable
  is equivalent to that      
  ${\rm dim}\,  {\rm Im } \gamma_g < \infty  $.

 The Heisenberg group associated to the infinite dimensional symplectic space $ F(( t))^{2N}$ is 
  \[ H= F(( t))^{2N} \times F \]
 with the group structure given by 
 \[  ( x_1 ,  k_1 ) ( x_2  ,  k_2 ) =  (x_1 + x_2 , \frac 1 2 \la x_1 , x_2 \ra + k_1 + k_2 ) . \] 
 The group $ \Sp_{2N} ( F((t))) $ acts on $H$ (from the right) by $ (x , k) \cdot g = ( x \cdot g  ,  k)$.

 We call  a complex valued function on $X_- $  a Schwartz function if its restriction to each finite dimensional subspace
   is a Schwartz function in the ordinary sense.  For example, if $F$ is a p-adic field with the ring of integers ${\cal O}$,
  the characteristic function of $ t^{-1} {\cal O} [ t^{-1}]^{2N} $ is a Schwartz function, and if $F$ is ${\Bbb R}$ or ${\Bbb C}$,
  $ q ( x ) $ is a ${\Bbb C}$-valued  quadratic form on $X_-$ with real part positive definite, then $ e^{- q ( x ) } $
  is a Schwartz function.  
Let 
  \[  {\cal S } ( X_-) \]
 denote the space of Schwartz functions on $X_-$.  We view $X_-, X_+$ and $F$ as  subgroups of $H$ by the embedding
    $v\in X_-\mapsto ( v, 0), v^*\in X_+ \mapsto ( v^* , 0 ) , k\in F \mapsto ( 0 , k )$.
The Heisenberg group $H$ acts on ${\cal S} (X_-)$ by the following:  for $k\in F$,  $v\in X_-$, $v^*\in X_+$,
 \begin{eqnarray}\label{hact} 
  (  k  \cdot f )  ( x ) &=& \psi ( k)  f ( x )   \\ 
  (   v  \cdot f ) ( x ) &=& f ( x + v )  \nonumber \\
  (  v^* \cdot  f ) ( x ) &=& \psi ( \la x , v^* \ra ) f ( x ) \nonumber 
\end{eqnarray}

We recall 
the results about the Weil representation for loop groups proved in \cite{Z}.

\begin{theorem}\label{act}  For each  $ g\in \Sp( F(( t))^{2N} , X_+ )$ with decomposition (\ref{gform}), 
   and a choice of Haar measure on ${\rm Im} \gamma $, we define an operator
 $T_g$ on ${\cal S}(X_- )$ by 
\begin{equation}\label{action}
( T_g f ) ( x ) = \int_{{\rm Im} \gamma}    S_{ g}( x + x^* )   f ( x \alpha + x^* \gamma ) d ( x^*\gamma  ), \end{equation}
where 
\[ S_{g} ( x + x^*) = \psi \left( \frac 12 \la x\alpha , x \beta \ra +
          \frac 12 \la x^* \gamma , x^* \delta \ra + \la x^* \gamma , x \beta \ra \right) ;\]
 $\psi $ here is a non-trivial additive character of $F$;
 then for each $h\in H$,
 \begin{equation}\label{2a}  T_g^{-1}  h T_g =  h\cdot g . \end{equation}
And $g\mapsto T_g$ gives a projective representation of $\Sp ( F(( t))^{2N}, X_+ )$ on ${\cal S} (X_- )$.
 \end{theorem}

 By restriction, we have a projective representation of $\Sp_{2N} ( F(( t)))$ on ${\cal S} ( X_- ) $.  
 The Steinberg symbol for the representation is given by the following: suppose $\sigma_1, \sigma_2 \in F[[t]]^*$
  with constant terms $c_1$ and $c_2$ respectively,  
  we have 
\begin{eqnarray}
 ( t^{2m}\sigma_1  ,  t^{2n} \sigma_2   )   &=&
     | C (t^{2m}\sigma_1 , t^{2n}\sigma_2 )|^{-\frac 12 }   \label{symbol}
 \\
  ( t^{2m }\sigma_1 ,  t^{2n+1} \sigma_2   )  &=&
    \frac { \gamma ( c_2 , \psi ) }{\gamma ( c_1 c_2 , \psi ) }  | C (t^{2m}\sigma_1 , t^{2n+1}\sigma_2 ) |^{-\frac 12 } \nonumber  \\
 ( t^{2m+1 }\sigma_1 ,  t^{2n} \sigma_2   )  &=&
    \frac { \gamma ( c_1 , \psi ) }{\gamma ( c_1 c_2 , \psi ) }  | C (t^{2m}\sigma_1 , t^{2n+1}\sigma_2 ) |^{-\frac 12 }  \nonumber \\
 ( t^{2m +1}\sigma_1 ,  t^{2n+1} \sigma_2   )  &=&
     { \gamma (c_1, \psi ) \gamma (c_2 , \psi ) } | C (t^{2m+1}\sigma_1 , t^{2n+1}\sigma_2 ) |^{-\frac 12 }  , \nonumber
\end{eqnarray}
  where $C(f_1 , f_2 )$ denotes the tame symbol of $f_1 $ and $f_2$ given by 
\[ C ( a , b ) =  (-1)^{v(a) v( b ) }  \frac {a^{v(b)}} { b^{v(a)}} |_{t=0}, \]
  and $\gamma ( c , \psi )$ denotes the 
  Weil index of $c\in F$ with respect to $\psi $ defined by the condition that the Fourier transform 
 \[  {\cal F} ( \psi ( \frac 12 c x^2 ) ) \]
 of the distribution $ \psi ( \frac 12 c x^2 )$ equals to 
 \[ \gamma ( c , \psi ) | c|^{-\frac 12 } \psi ( - c^{-1} x^2 ), \]
   see \cite{W1}. 

\

 We fix a Borel subgroup $B_0$ of ${\rm Sp}_{2N} $ as the stabilizer of the flag 
 \[ {\rm span}\{e_1\} \subset \dots \subset {\rm span} \{ e_1 , \dots , e_N\},\] 
 where $e_i$ denotes the vector with $i$-th coordinate $1$ and other coordinates $0$.
 We fix the maximal torus $A_0$ of $Sp_{2N} ( F ) $ consisting of diagonal elements, and let $\Delta_0$ be the corresponding set of roots of ${\rm Sp}_{2N}$.
  We  fix a Chevalley basis of the Lie algebra of ${\rm Sp}_{2N}$.  The group 
    $Sp_{2N} (F((t)))$ can be described as the group  
generated by root vectors $x_{\alpha } ( a ) $, where $\alpha \in \Delta_0 $ and $a\in F((t))$, and 
    the relations (\ref{llg1}) (\ref{llg2}) and (\ref{llg3}) below:
 \begin{equation}\label{llg1}
      x_{\alpha } ( a_1 )  x_{\alpha } ( a_2 )  = x_{\alpha } ( a_1 + a_2  ). 
\end{equation}
    If $\alpha$ and $\beta $ are roots and $\alpha +\beta \ne 0$, then
  \begin{equation}\label{llg2}
           x_{\alpha} (a ) x_{\beta}(b) x_{\alpha} (a )^{-1} x_{\beta}(b)^{-1} =
     \Pi x_{i\alpha +j\beta } ( c_{ij} a^i b^j),
  \end{equation}
 where the product is over all the roots $i \alpha + j \beta $, $i>0, j> 0$
 and  the coefficients $c_{ij}\in {\bf Z}$  are given  in terms of the Chevalley basis of $\frak g$.
 If $\alpha + \beta $ is not a root, then the right hand side is $1$.  
 See \cite{S} for the precise meaning of the right hand side. 
    For $a \in F((t))^*$, we set  
 \[ w_{\alpha }(a )= x_{\alpha}(a  )x_{-\alpha}(-a^{-1})x_{\alpha}(a ) \]
  and 
   \[ h_{\alpha} ( a) = w_{\alpha }(a )w_{\alpha}(1)^{-1} ;  \]
 then 
    \begin{equation}\label{llg3} 
   h_{\alpha} ( a_1) h_{\alpha} ( a_2 ) = h_{\alpha } ( a_1 a_2 ) .
  \end{equation}
 When $G=SL_2$, there are two roots $\alpha $ and $-\alpha $,  the relations (\ref{llg2}) above is replaced by 
 \begin{equation}\label{llg4}     w_{\alpha} ( a )  x_{\alpha }(b) w_{\alpha } ( -a ) =  x_{- \alpha }(-a^{-2}b). \end{equation}
Let $\widetilde{\Sp}_{2N} (( F((t)))$ denote the Steinberg group defined by the symbol (\ref{symbol}).  It is
  generated by root vectors $x_{\alpha } ( a ) $ ($\alpha \in \Delta_0, a \in F(( t))$) and ${\Bbb C}^*$  with relations 
 (\ref{llg1}) (\ref{llg2}) and 
   \begin{equation}\label{llg5}
    h_{\theta }  ( a ) h_{\theta } ( b ) h_{\theta } ( a b )^{-1} = (a , b),
 \end{equation}  
 where $(a , b) $ is the symbol in (\ref{symbol}), and $\theta $ is the longest root of ${\rm Sp}_{2N}$,   
and ${\Bbb C}^*$ is in the center.   
   
 In the case of  ${\rm SL}_2 $, $\widetilde{\rm SL}_2 ( F ((t)))$ is generated by $x_{\alpha } ( a ), x_{-\alpha } ( a )$ and
 ${\Bbb C}^*$ with relations
  (\ref{llg1}), (\ref{llg4}), (\ref{llg5}), and ${\Bbb C}^*$ is in the center.
We call $\widetilde{\Sp}_{2N} (( F((t)))$ the metaplectic loop group of ${\Sp}_{2N} (( F((t)))$.
  The group $\widetilde{\Sp}_{2N} (( F((t)))$ is a central extension of ${\Sp}_{2N} (( F((t)))$:
\[     1 \to {\Bbb C}^* \to \widetilde{\Sp}_{2N} ( F((t))) \to {\Sp}_{2N} ( F((t))) \to 1 .\]
 We denote the image of $g\in \widetilde{\Sp}_{2N} (( F((t)))$ in ${\Sp}_{2N} (( F((t)))$ by $\bar g $.  
  Since the symbol (\ref{symbol}) is trivial on the subgroup $ F^* \times F^* \subset F(( t))^* \times F((t))^*$,
    the elements $x_{\alpha } ( a ) $ ($a\in F$) generate a  
  subgroup of  $\widetilde{\Sp}_{2N} (( F((t)))$  isomorphic to ${\rm Sp}_{2N} ( F )$, so  we will regard $\Sp_{2N} ( F)$ as a 
subgroup of  $\widetilde{\Sp}_{2N} (( F((t)))$.

 To describe the action  $\pi :   \widetilde{\Sp}_{2N} (( F((t))) \to GL ( {\cal S} ( X_- ) ) $, 
   it is sufficient to describe the action of the generators $x_{\alpha } ( a )$ and of $c\in {\Bbb C}^*$.
  If $a\in F[[ t]]$, then $ g \doteq \bar x_{\alpha } ( a ) \in Sp_{2N} ( F[[t]] )$,
  ${\rm Im } \gamma_g =\{ 0 \}$,   $ \pi ( x_{\alpha } (a ) ) $ is $T_g$ in
  (\ref{action}) with the Haar measure on ${\rm Im \gamma}_g $ as the counting measure, i.e., the volume of $\{ 0 \}$ is $1$. 
   For general $a\in F((t))$, we can find a diagonal 
    element of type 
 \[  t^k \doteq {\rm diag} ( t^{k_1} , \dots , t^{k_N} , t^{-k_1} , \dots , t^{-k_N} ) \]
  such that   $ t^k  \bar x_{\alpha } ( a )  (t^k)^{-1} = \bar x_{\alpha } ( t^{( \alpha , k )} a ) 
   \in {\rm Sp}_{2N} ( F[[ t ]] ) $; then we set 
  \begin{equation}\label{2a1}  \pi ( x_{ \alpha } ( a ) ) =  (T_{t^k})^{-1} \pi ( x_{\alpha } ( t^{( \alpha , k )}a ) ) T_{t^k} . \end{equation} 
 Note that $T_{t^k}$ is as in (\ref{action}) and the conjugation above is independent 
  of the choice of the Haar measure for ${\rm Im} \gamma_{t^k} $.
 And  $c \in {\Bbb C}^*$ acts as scalar multiplication by $c$.

\

\begin{lemma}\label{lemma2.1}\label{lemma2.1a} The representation of $ \widetilde{\Sp}_{2N} ( F((t)))$ on ${\cal S} ( X_- ) $ is 
  faithful, i.e.,  $\pi :   \widetilde{\Sp}_{2N} ( F((t))) \to {\rm GL} ( {\cal S} ( X_- ) ) $
 is injective.
 \end{lemma}

\noindent   {\it Proof.}  By Theorem \ref{act}, we have,  for $g\in  \widetilde{Sp}_{2N} (( F((t)))$ and $h\in H$,
 \begin{equation}\label{compatible}
  \pi ( g )^{-1} h \pi ( g ) =        h \cdot {\bar g} . \end{equation}
  If $\pi ( g ) = 1$, then $h \cdot {\bar g} = h $ for all $h\in H$,    
so $\bar g = 1 $.   This means $ g \in {\Bbb C}^*$, so $g=1$. 
\hfill $\Box $

\

  When we need to indicate the dependence of the Weil representation $\pi $ on the additive character $\psi$, we
  write $\pi$ as $\pi_{\psi }$.   If $\psi '$ is related to $\psi $ by the relation
  $\psi' ( x ) = \psi ( b^2  x )$, then  we have 

\begin{lemma}\label{2.3}    The map $f (x ) \mapsto f(b x )$ is an isomorphism from Weil representation
  $\pi_{\psi}$ to $\pi_{\psi'}$. 
\end{lemma}

\noindent   {\it Proof.}  We denote the map $f ( x) \mapsto f ( b x ) $ by $\Phi$. It is direct to 
 check that 
 \begin{equation}\label{2f}
   \pi_{\psi'} ( x_{\alpha } (a )  )   \Phi  =  \Phi  \pi_{\psi} ( x_{\alpha } ( a )  )
 \end{equation} 
for  $ a \in F[[t]]$.    We then check (\ref{2f}) holds for  $ a \in F ((t ))$ using (\ref{2a1}).
\hfill $\Box $.

\

The reprametrization group of $F(( t))$ is by definition
\[{\rm Aut} F(( t)) = \{  \sum_{i=1}^{\infty} a_i t^i \in F[[t]] t \, | \, a_1\ne 0 \} ,\] 
with the group operation $  (\sigma_1 * \sigma_2 )(t) = \sigma_2 ( \sigma_1 ( t ) ) .$
It acts on $F((t))$ from the right  by 
 \[  a ( t ) \cdot \sigma (t)   = a ( \sigma^{-1} ( t ) ) . \]
 And it acts on the space $F((t))dt$ of formal 
 $1$-forms (from the right)  by 
\[  a ( t ) dt  \cdot \sigma = a ( \sigma^{-1} ( t ) ) \sigma^{-1} ( t)  '  dt . \]
We  view the first $n$ components in $ X= F((t))^{2n}$ as elements in $F(( t))$ and
 the last $n$ components as elements in $F(( t)) dt $ (without writing $dt$), then
  $ {\rm Aut} F(( t))$ acts on $ F((t))^{2n} $ by 
  \begin{eqnarray}
  &\,&  ( a_1 ( t ), \dots , a_n ( t) , a_{n+1}(t) , \dots , a_{2n}(t) )  \cdot \sigma   \nonumber \\
  &=& (a_1 ( \sigma^{-1} ( t ) ), \dots , a_n ( \sigma^{-1} ( t ) ), a_{n+1 }( \sigma^{-1} ( t ) ) \sigma^{-1} ( t )' ,  \dots ).
 \nonumber
   \end{eqnarray}
Since the residue of an $1$-form is independent of the local parameter,  the action preserves the symplectic form. 
 We have an embedding $ {\rm Aut} F(( t)) \subset  \Sp ( F(( t))^{2N} )$. And when we write  
    $\sigma  \in {\rm Aut} F((t))$ as in (\ref{gform}),  it is clear that  $\gamma = 0 $,  
  so we may view  ${\rm Aut} F((t))$ as a subgroup of $ \Sp ( F(( t))^{2N} , X_+ )$.
Using (\ref{action}),  $\sigma $ acts on ${\cal S} (X_- )$ by  
\[  ( \pi ({\sigma}) f ) (x ) = \psi ( \frac 1 2 \la x \alpha , x \beta \ra ) f ( x \alpha ) . \] 
 It is easy to check that  
 \[  \pi ({\sigma_1 }) \pi ({\sigma_2 } )  = \pi ({\sigma_1 * \sigma_2 } ).\]
The multiplicative group $ F^*$ is a subgroup of $ {\rm Aut} F(( t))$ by the embedding $ c \mapsto ct $, the 
 action of $F^*$ on ${\cal S} (X_-)$ is given by 
\begin{eqnarray}\label{deop}
&\, &  (c \cdot f )( x_1 ( t ), \dots , x_n ( t) , x_{n+1}(t) , \dots , x_{2n}(t) )  \\
 &\, & =  f ( x_1 ( c^{-1} t ) , \dots , x_n ( c^{-1} t  ), c^{-1} x_{n+1 }( c^{-1} t  ),  \dots , c^{-1} x_{2n}( c^{-1} t) ). \nonumber
\end{eqnarray}

The group
  ${\rm Aut} F((t))$ also acts on $\Sp_{2N} ( F((t)))$ as automorphisms in the following way:
  for $\sigma (t)\in {\rm Aut} F((t))$, $g\in \Sp_{2N} (( F((t)))$, we write
   \[  g = \left( \begin{array}{cc}  a  &  b \\ c & d  \end{array} \right) \]
  according to the decomposition $ F(( t))^{2N} = F (( t))^{N} \oplus  F((t))^{N}$, where
  as above, the entries in the first component are functions and the entries of the second component are  $1$-forms, 
  so the blocks $a$ and $d$ are  $N\times N$-matrices with entries as functions, the block $b$ is an
  $N\times N$-matrix with entries as $1$-forms, and the block $c$ is an $N\times N$-matrix with entries
  as vector fields;  $\sigma (t)$ changes an entry $k(t)$ in block $a$ or $d$ to $ k ( \sigma (t))$,  changes
  an entry $k(t)$ in block $b$ to $k(\sigma (t)) \sigma(t)' $, and changes an entry $k(t)$ in block $c$ to 
  $k(\sigma (t))(\sigma  (t)')^{-1}$.  This action is compatible with the action on $H$.  
   Since the ${\rm Aut} F((t))$-action on $F(( t))$ preserves the symbol in \ref{symbol}, 
 the  ${\rm Aut} F((t))$-action on  $\Sp_{2N} (( F((t)))$
 lifts to an action
   on   $\widetilde{\Sp}_{2N} ( F((t)))$ and therefore an action of 
  ${\rm Aut} F((t))$ on the semi-direct product  $\widetilde{\Sp}_{2N} ( F((t))) \rtimes H $. And ${\cal S} (X_- )$
  is a representation of the semi-direct product group
  ${\rm Aut} F((t))\rtimes (\widehat{\Sp}_{2N} ( F((t))) \rtimes H  )$. 

\

 \noindent { \bf  2.2. Bruhat and Iwasawa decompositions . } 
  
 Let $B$ be the subgroup of $\widetilde {Sp}_{2N} ( F((t)) )$ which consists of elements $g$ such 
  that $ {\bar g} \in  \Sp_{2N} ( F[[t]] )$ and $ {\bar g} \, {\rm mod} \, t $ is in $B_0$.  We call $B$ a Borel 
  subgroup of $\widetilde {\Sp}_{2N} ( F((t)) )$. It is clear the center ${\Bbb C}^* \subset B$.
  Let $N$ be the subgroup generated by $w_{\alpha } ( a )$ with 
 $\alpha \in \Delta_0$ and $a \in F((t))^*$.   Then $(B, N)$ is a BN-pair  for   $ \widetilde{\Sp}_{2N} ( F((t))) $ with the affine 
 Weyl group
  ${\widehat W}$ as the Weyl group. This can be proved using the pull-back of the 
 standard BN-pair for $ {\Sp}_{2N} ( F((t)) )$ under the map   $\widetilde {\Sp}_{2N} ( F((t)) ) \to {\Sp}_{2N} ( F((t)) )$.
Recall that $\widehat W$ is the semi-direct product of the 
  Weyl group $W$ and the coroot lattice $Q^{\vee} $ of $\Sp_{2N}$.  
  We have the Bruhat decomposition for $ \widetilde{\Sp}_{2N} ( F((t))) $:
\[      \widetilde{\Sp}_{2N} ( F((t)))  = \sqcup_{w} B w B  ,\]  
   where $w$ runs through all elements in $\widehat W$.

  We wish to define a ``maximal compact subgroup'' $K$ for $ \widetilde {\Sp}_{2N} ( F((t))) $. If $F$ is a p-adic field with ring of integers
 ${\cal O}$,  we let $K$ be the subgroup generated by $x_{\alpha } ( a )$ with $\alpha \in \Delta_0$ and $a\in {\cal O} (( t))$. 

 For each affine real root $ n \delta + \alpha $ ($n\in {\Bbb Z}$, $\alpha \in \Delta_0$),  we call the 
 group defined by
   \[  \{  x_{\alpha } ( c t^n ) \, | \, c \in F \} \]
 the root group for $n \delta + \alpha $. In particular, 
for  $\alpha_0 = \delta - \theta $,   the extra simple root of affine ${\Sp}_{2N}$,
   we have  the root groups 
  \[  \{  x_{-\theta } ( c t )  \, |  \, c \in F  \}, \, \, \, \, \, \{  x_{\theta } ( c t^{-1} )  \, |  \, c \in F  \}  \, \, \, \, \,    
    \]
 for $\alpha_0$ and $-\alpha_0$.   We write $L_+ ( c ) = x_{-\theta } ( c t )$ and 
  $L_- ( c ) = x_{\theta } ( c t^{-1} )$, and put for $b\in F^*$, 
\[  w ( b ) = L_+ ( b ) L_-( -b^{-1} ) L_+ ( b ) , \, \, \, \,  h ( b ) = w( b ) w( 1 ) ^{-1}. \] 
A direct calculation using the generating relations of
 $\widetilde{\Sp}_{2N} ( F ( (t)))$ gives the following:
  \[  L_{\pm} ( c_1 )   L_{\pm } ( c_1 ) = L_{\pm } ( c_1 +c_2  ) ,\] 
 \[  w ( b )  L_+ ( c ) w ( - b ) = L_- ( - b^{-2} c ) , \]
\[   h( b_1  ) h ( b_2 ) h ( b_1 b_2 ) ^{-1} = \frac { \gamma ( b_1 , \psi ) } {   \gamma ( 1 , \psi ) }
   \frac { \gamma ( b_2 , \psi ) } {   \gamma ( 1 , \psi ) }  \frac { \gamma ( 1, \psi ) } {   \gamma ( b_1 b_2  , \psi ) }.\]
Notice that the right hand side of the last identity is the Hilbert symbol for $b_1$ and $b_2$ (see \cite{W1} page 176) .
  So the subgroup $G_{\alpha_0} $  generated by $x_{ -\theta } ( c t)$ and $x_{\theta } ( c t)$  ($c\in F$) 
   is a  central extension of ${\rm SL}_2 (F)$ under the map
\[    x_{\theta } ( c t^{-1} ) \mapsto   \left(
 \begin{array}{cc}  1  &  0 \\ c &  1  \end{array} \right), \, \, \, \,  
x_{-\theta } ( c t ) \mapsto   \left(
 \begin{array}{cc}  1  &  c \\ 0 &  1  \end{array} \right). \]
 When $F={\Bbb C}$, this map is an isomorphism; and when $F\ne {\Bbb C}$, the kernel 
  is $\{ \pm 1 \}$,  and $G_{\alpha_0} $ is a two-fold cover of ${\rm SL}_2(F)$ given 
  by the Hilbert symbol of $F$, which is called the metaplectic group for $SL_2 ( F)$.  
  We let $K_{\alpha_0}$ 
denote the standard, maximal compact subgroup of $G_{\alpha_0} $. For example,
 if $F={\Bbb C}$, $K_{\alpha_0}$ is $SU_2$.  For each simple root $\alpha_i $ ( $i=1 , \dots , N $)
 we let $ K_{\alpha_i } $ be the standard maximal compact subgroup for the $SL_2 ( F_v )$   
  associated to $\alpha_i $.

  For $F = {\Bbb R}$ or ${\Bbb C}$, we first let $K_{\rm fin} = {\rm SO}_{2N} \cap \Sp_{2N} ( {\Bbb  R} )$ or  
$K_{\rm fin}= {\rm SU}_{2N} \cap  \Sp_{2N} ( {\Bbb C}) $ according to whether $F = {\Bbb R}$ or ${\Bbb C}$.   
  Let $K $ be the subgroup generated by $K_{\rm fin}$ and $K_{\alpha_0}$.  So we have chosen a ``maximal compact subgroup" 
  $K$ of $\widetilde{\Sp}_{2N} ( F((t)))$ for each local field $F$ of characteristic $0$.

  For a p-adic field $F$, let $K'$ be the subgroup generated by $K_{\rm fin}= \Sp_{2N} ({\cal O})$
  and $K_{\alpha_0}$, it is easy to prove that $ K' \subset K$.

  By use of the BN-pairs mentioned above and the method as in \cite{S},  we can prove  the Iwasawa decomposition

\begin{equation}\label{I1}    \widetilde{\Sp}_{2N} ( F((t)))  =   B  K  .  \end{equation}
 
\noindent  For a p-adic field $F$, we also have 
  
\begin{equation}\label{I2}    \widetilde{\Sp}_{2N} ( F((t)))  =   B  K'  .  \end{equation}

\begin{lemma}\label{Lemma2.4a}  There is a splitting homomorphism $Sp_{2N} ( F[[ t]] ) \to  \widetilde {Sp}_{2N} ( F((t))) $.
\end{lemma}

\noindent {\it  Proof.} Consider the Weil representation of $\widetilde{\Sp}_{2N} ( F((t))) $.  For
   each $g \in \Sp_{2N} ( F[[ t]] )$, we have ${\rm Im} \, \gamma_ g =\{ 0 \}$, we taking the counting measure 
  in the formula (\ref{action}) for $T_g$,  then it is easy to see that $T_{g_1}T_{g_2} = T_{g_1 g_2}$
  for $g_1, g_2 \in  \Sp_{2N} ( F[[ t]] )$.  The map $g \mapsto T_g$ defines the desired splitting 
$Sp_{2N} ( F[[ t]] ) \to  \widetilde{Sp}_{2N} ( F((t))) $.   \hfill $\Box $

\

With this lemma, we shall regard $Sp_{2N} ( F[[t]])$ as a subgroup of  $\widetilde{Sp}_{2N} ( F((t))) $.
Also it is clear that  $B$ is  a subgroup of  $ {\Bbb C}^* Sp_{2N} ( F[[t]])$.

\begin{lemma}\label{phi0RC}   If $F= {\Bbb  R}$ or $ {\Bbb C}$,  there is a non-zero element
  $\phi_{0 } \in  {\cal S} ( t^{-1} F[t^{-1}]^{2N})$ fixed 
  by $K$ up to a scalar.
\end{lemma}

\noindent {\it Proof.} If $F= {\Bbb R}$, by Lemma \ref{2.3}, we may assume that $\psi ( x )
  =e^{\pm 2 \pi i x } $. The result then follows from  \cite{Z}, Section 4, 
  where
   $\phi_0 = {\rm exp} ( \pi i (x, x \Omega ))$
  for $\Omega = i I $ (see, Section 4,  \cite{Z} for the definition of $I$).   If $F={\Bbb C}$, we may assume that
 $\psi ( z ) =e^{ 2\pi i ({\rm tr} z ) } = e^{ 4 \pi i ({\rm re} \, z )  }$.     We view 
  $ \Sp_{2N} ( {\Bbb C} (( t)))$ as a subgroup of $ \Sp_{4N} ( {\Bbb R} (( t))) $ under 
   then identification $ {\Bbb C} = {\Bbb R}^2 $, 
    $z = x+ i y \mapsto (x , y ) \in {\Bbb R}^2 $.   
 We may view  
   $\widetilde {\Sp}_{2N} ( {\Bbb C} (( t)))$
 as a subgroup of $\widetilde{Sp}_{4N} ( {\Bbb R} (( t))) $.
  Then the ``maximal compact subgroup" of $\widetilde {\Sp}_{2N} ( {\Bbb C} (( t)))$ is a subgroup of the maximal subgroup of  
 $\widetilde{\Sp}_{4N} ( {\Bbb R} (( t))) $.
   The Weil representation for  $ \widetilde{\Sp}_{4N} ( {\Bbb R} (( t)))$
    restricts to the Weil representation of  $ \widetilde{\Sp}_{2N} ( {\Bbb C} (( t)))$.
  The function $\phi_0$ is fixed by the ``maximal compact" subgroup of 
     $ \widetilde{\Sp}_{2N} ( {\Bbb C} (( t)))$ up to a scalar.   \hfill $\Box $

\

We expect that $\phi_0$ in the Lemma is actually fixed by $K$.

\

\begin{lemma}\label{phi0padic} If  $F$ is a 
  non-Archimedean local field with residual characteristic not equal to $2$ and ring of integers  ${\cal O}$,
  and if the conductor of $\psi $ is ${\cal O}$, i.e., 
 \[ {\cal O} = \{ x \in F \, | \, \psi ( x b) =1,  \, \, \, {\rm for \, all} \, \, b \in {\cal O} \}. \]
Then the characteristic function $\phi_0 $
 of $ t^{-1} {\cal O} [t^{-1}]^{2N} $ is fixed by $K$.
\end{lemma}  

\noindent {\it Proof.}  It is enough to check $x_{\alpha } ( a ) $ for $a \in {\cal O} ((  t))$ fixes $\phi_0 $.
 If $a \in {\cal O} [[t]]$, it is easy to check  $x_{\alpha } ( a ) $ fixes $\phi_0 $. 
 For $ a \in {\cal O} ((t))$, we use formula (\ref{2a1}).  Since $T_{t^k}$ is a partial Fourier transform, it fixes $\phi_0 $,
 and $\pi ( x_{\alpha } ( t^{ ( \alpha , k ) } a )) $ fixes $\phi_0 $ as $t^{ ( \alpha , k ) } a\in {\cal O} [[t]] $, so 
  $x_{\alpha } ( a ) $ fixes $\phi_0 $.
\hfill $\Box $

\

We shall fix $\phi_0 \in {\cal S} ( t^{-1} F [ t^{-1} ] ) $ as described in Lemma \ref{phi0RC} and Lemma \ref{phi0padic}. 
 Lemma \ref{phi0padic} implies  

\begin{lemma}  If  $F$ is a 
  non-Archimedean local field with residual characteristic not equal to $2$ and the conductor of $\psi $ is 
  ${\cal O}$, the  map $\widetilde{\Sp}_{2N}  ( F (( t)))\to {\Sp}_{2N} ( F (( t)))$ maps $ K$ to
  $\Sp_{2N} ( {\cal O} (( t)))$ isomorphically.
\end{lemma}

{\it Proof.}  It is clear that the map is surjective onto  $Sp_{2N} ( {\cal O} (( t)))$.
  Since the kernel is in ${\Bbb C^*}\cap K$, and $K$ fixes $\phi_0$, we must have 
 ${\Bbb C}^* \cap K= 1 $. \hfill $\Box$.

\

\noindent { \bf  2.3. Adelic groups and representations.}  In this section, we assume $F$ is a number 
  field. By Section 2.1, we have, for each place $v$ of $F$, a representation 
   of the semi-direct product $ {\rm Aut } ( F_v ((t))) \rtimes  \widetilde{\Sp}_{2N} ( {F}_v (( t))) $
   on the space $ {\cal S} ( X_{-, v } )$, where  $ X_{-, v } = t^{-1}{F_v } [ t^{-1} ]^{2N} $.
   In this section we define the adelic metaplectic loop group for ${\Sp}_{2N}$ and define its Weil representation.

  Let ${\bf A}$ denote the ring of adeles of $F$, and $\psi =\Pi \psi_v$ be a 
  non-trivial character of ${\bf A} / F$.  
 For a non-Archimedean place $v$,
  ${\cal O}_v$ denotes the ring of integers.  We let
\[  
               {\bf A} \la t \ra  = \{ (a_v) \in \Pi_v F_v (( t)) \, | \,  a_v \in {\cal O}_v (( t)) \, \, {\rm for \, \, almost \, \,  all } \, v \},
 \]
  
\[  {\bf A} \la t \ra_+  =  {\bf A} \la t \ra  \cap {\bf A} [[ t]] ,\]
 \[  {\bf A} \la t \ra_-  =  \{ (a_v) \in \Pi_v t^{-1} F_v [t^{-1} ] \, | \,   a_v \in {\cal O}_v (( t)) \, \, {\rm for \, \,
almost \, \,  all } \, v \}.\]
It is clear that 
   \[  {\bf A} \la t \ra = {\bf A} \la t \ra_+ \oplus  {\bf A} \la t \ra_- , \, \, \,\,  t^{-1} {\bf A} [t^{-1} ] \subset  {\bf A} \la t \ra_- .\]
And we let  
 \[   F \la t \ra   =  {\bf A} \la t \ra  \cap  F(( t)) .\]
  An element $a \in F((t))$ is in $F\la t \ra $ iff for almost all finite places $v$, $a \in {\cal O}_v (( t)) $.

\begin{lemma}  $F \la t \ra $ is subfield of $F(( t))$.
\end{lemma}

\noindent {\it Proof.} It is enough to prove that if $ a= \sum_{ i\geq l} k_i t^i \in F \la t \ra -\{ 0 \}$, then 
   $a^{-1} \in  F \la t \ra$.   Let $S$ be the  set of finite places $v$ such that $ a \notin {\cal O}_v (( t)) $.
  We may assume $k_l \ne 0 $. Let $S'$ be the finite set of finite places $v$ such that $ a_l $ is not a  unit of ${\cal O}_v$.
 It is clear that both $S$ and $S'$ are finite sets.  
Then $a^{-1} \in {\cal O}_v (( t)) $ for $ v \notin S \cup S' $. 
 \hfill $\Box $

\

   The adelic loop group  for $\Sp_{2N}$ without central extension  is defined as
  $ \Sp_{2N} (  {\bf A} \la t \ra  )$.  It is clear that 
  $ \Sp_{2N} (  {\bf A} \la t \ra  )$ is the restricted product of $  \Sp_{2N} (  F_v (( t)) )$
 with respect to ``maximal compact" subgroups $ \Sp_{2N} (  {\cal O}_v   (( t ))  )$. 
  The adelic metaplectic loop group $  \widetilde{\Sp}_{2N} ( {\bf A } \la t \ra ) $ is defined as the 
  restricted product $ \Pi'_v   \widetilde{\Sp}_{2N} ( {F}_v (( t)))$ with respect to ``maximal compact" subgroups
  $K_v$ defined in Section 2.2. Clearly we have the exact sequence

   \begin{equation}\label{2.3.1}    1\to \oplus_v {\Bbb C}^* \to    \widetilde {\Sp}_{2N} ( {\bf A} \la t \ra ) 
            \to   {\Sp}_{2N} ( {\bf A} \la  t \ra ) \to 1 . \end{equation}

 We  can also define  the adelic  group for ${\rm Aut} F((t))$. 
   For ${F}_v ={\Bbb R}$,   we define the 
  maximal ``compact" subgroup of ${\rm Aut} ( {\Bbb R} ((t)))$ as $\{ -t , t \}$.
   If ${ F}_v ={\Bbb C}$,   we define the 
  maximal ``compact" subgroup of ${\rm Aut} ( {\Bbb C} ((t)))$ as $\{ c t \, | \, |c|=1  \}$.
   If  $F_v $ is a non-Archmedean local field,  
     we define the 
  maximal ``compact" subgroup of ${\rm Aut} ( {F}_v ((t)))$ as the subgroup 
   consisting of elements
  $\sum_{i=1}^{\infty} c_i t^i $ with $c_1 \in {\cal O}_v^*$ and $c_i\in {\cal O}_v$ for all $i\geq 2$.
 It is easy to check that the above maximal ``compact" subgroups of $ {\rm Aut} ( { F}_v ((t)))  $ 
  preserves $K_v$ of   $\widetilde{Sp}_{2N} ( {F}_v (( t)))$ defined in Section 2.2.   
The adele group for ${\rm Aut} F((t))$ is 
\[ {\rm Aut } {\bf A}\la t \ra  =  \Pi'_v   {\rm Aut } ( {F}_v ((t))), \]
  where the restricted product is with respect to the ``compact" subgroup of ${\rm Aut } ( {F}_v ((t)))$ defined as above.
  It is clear that ${\rm Aut } {\bf A} \la t \ra $ acts on   $\widetilde {Sp}_{2N} ( {\bf A} \la t \ra  ) $, so we have 
 semi-direct product ${\rm Aut } {\bf A} \la t \ra  \rtimes \widetilde {Sp}_{2N} ( {\bf A} \la t \ra  )$.

  For each place $v$, as in Section 2.1, we have the Weil representation $ {\cal S} ( X_{ - , v } ) =  {\cal S } (  t^{-1} F_v [ t^{-1} ]^{2N} )$ 
of $ \widetilde{\Sp}_{2N} ( {F}_v (( t)))$. 
Since for almost all places $v$, there is $\phi_{v, 0}$ (Lemma \ref{phi0padic}) fixed by
the local ``maximal compact" group $K_v$,  the restricted tensor product $ \otimes' {\cal S} ( ( F_v [t^{-1}] t^{-1} )^{2N} )   $
  with respect to $\{ \phi_{v, 0} \}$ is a representation of 
$ \widetilde {Sp}_{2N} ( {\bf A} \la t \ra  )$. 
  We call this representation the adelic Weil representation.
And note that for almost all places $v$, the maximal compact subgroup
  of ${\rm Aut} F_v (( t))$ fixes $\phi_{v, 0 }$, so  ${\rm Aut } {\bf A} \la t \ra  \rtimes \widetilde {Sp}_{2N} ( {\bf A} \la t \ra  )$
  acts on $\otimes' {\cal S} ( ( F_v [t^{-1}] t^{-1} )^{2N} ) $,  it is clear that 
$\otimes' {\cal S} ( ( F_v [t^{-1}] t^{-1} )^{2N} ) $ can be regarded as a function space on   
  $  ({\bf A} [ t^{-1} ] t^{-1} )^{2N}$.

 \begin{lemma}   There is a splitting homomorphism $\Sp_{2N} ( F \la t \ra  ) \to    \widetilde {\Sp}_{2N} ( {\bf A} \la t \ra  ) $ 
  \end{lemma}

\noindent{\it Proof.}   Since  $F \la t \ra $ is a field,   $\Sp_{2N} ( F \la t \ra  ) $  is isomorphic to the group
 generated by the root subgroups $y_{\alpha } ( a ) $ ($\alpha \in \Delta_0$ , $a\in  F \la t \ra  $) with the standard relations
  (\ref{llg1}) (\ref{llg2}) and (\ref{llg3}) ( (\ref{llg4})if $N=1$). 
  On the other hand,  for $\alpha \in \Delta_0, a\in F\la t \ra $, we have
  $x_{\alpha } ( a ) \in \widetilde {\Sp}_{2N} ( {\bf A} \la t \ra  )$. Because of  
  the product formula $ \Pi_v |c|_v = 1 $  and because of the product formula for the Weil index  $\Pi_v \gamma ( c , \psi_v) = 1 $ when 
  $c\in F^*$,  so we have, for $f_1 , f_2 \in F\la t \ra $,  the product formula  
\[ \Pi_v  ( f_1 , f_2  )_v = 1 ,\]  
 where $ (f_1, f_2)_v$ denotes the symbol (\ref{symbol}) for $F_v((t))$. 
So  $x_{\alpha } ( a )$ ($a\in F\la t \ra$)  satisfies the standard relations  (\ref{llg1}) (\ref{llg2}) and (\ref{llg3})
 ( (\ref{llg4})if $N=1$) . Therefore 
\[ y_{\alpha } ( a ) \mapsto x_{\alpha } ( a ) \]
 is the desired splitting.    \hfill $\Box$

\

 From now on, we regard  $\Sp_{2N} ( F \la t \ra  )$ as a subgroup of $ \widetilde {\Sp}_{2N} ( {\bf A} \la t \ra  ) $
 whenever convenient.
    
   We also have  the Iwasawa decomposition for adelic groups 
\begin{equation}\label{ad}  \widetilde {\Sp}_{2N} ( {\bf A} \la t \ra  )       = B_{\bf A}  K_{\bf A}, \end{equation}
 where   $B_{\bf A}$ is the restricted product of groups $B_v$ with respect to $B_v\cap K_v $
  and $K_{\bf A} = \Pi_v K_v $.

\

\section{ Theta Functional }

\noindent  We continue to assume $F$ is a number field.  In this section, we first 
  introduce a certain function space   ${\cal E } ( t^{-1} {\bf A} [ t^{-1} ]^{2N} )$ on 
 $ t^{-1} {\bf A} [ t^{-1} ]^{2N} $ that is closed under the action of  
 the adelic metapletic group $\widetilde{\rm Sp}_{2N} ( {\bf A} \la t \ra ) $.  For  $T \in {\rm Aut} ( {\bf A} \la t \ra )$ satisfying 
  certain conditions,
 we  construct 
  theta functional $\theta :   T  {\cal E } ( t^{-1} {\bf A} [ t^{-1} ]^{2N} ) \to {\Bbb C}$ that is invariant
  under the action of ${\rm Sp}_{2N} ( F \la t \ra ) $.
  
  The definition of functional $\theta $ is similar to the classical case: 
  \[ \theta ( f ) =  \sum_{ r \in t^{-1} F [ t^{-1} ]^{2N} } f ( r ) .\]
  It is easy to see that the above summation is not convergent for arbitrary 
  $ f\in  \otimes' {\cal S} ( ( F_v [t^{-1}] t^{-1} )^{2N} ) $. Our main result (Theorem \ref{converge}) is
  that $\theta (f ) $ converges  for  
  $ f \in T {\cal E } ( t^{-1} {\bf A} [ t^{-1} ]^{2N} )$.

    \
  
   For a finite place $v$, a subgroup of the Heisenberg group $H_v = F_v (( t))^{2N} \times F_v $ is called a congruence subgroup
   if it contains   $\pi_v^k  {\cal O}_v (( t))^{2N}$ for some positive integer $k$, where $\pi_v \in {\cal O}_v$
  is a local prime.  

\begin{lemma}\label{lemma3.0}  If $v$ is a finite place, suppose $ \phi \in {\cal S} ( t^{-1} F_v [t^{-1} ]^{2N})$
  is fixed by  $\pi_v^k  {\cal O}_v (( t))^{2N}$. Then $\phi $ is invariant under the translation by elements in 
 $\pi_v^k   t^{-1}{\cal O}_v [ t^{-1} ]^{2N}$ and is  
supported in $ \pi_v ^{-k-l } t^{-1} {\cal O}_v [ t^{-1}  ]^{2N} $,
 where $\pi_v^{-l}{\cal O}_v $ is the conductor of $\psi_v $, 
\end{lemma}

\noindent{\it Proof.}  Since elements in $ t^{-1}{\cal O}_v [ t^{-1} ]^{2N}$ act on 
  ${\cal S} ( t^{-1} F_v [t^{-1} ]^{2N})$ by  translations (see (\ref{hact})), 
 $\phi $ is fixed by $\pi_v^k  t^{-1} {\cal O}_v [t^{-1} ] ^{2N}$ means precisely the first 
 claim. For the second claim,  since $\phi $ is fixed by every $ v^* \in \pi_v^k {\cal O}_v [[ t]]^{2N}$,
  we have 
  \begin{equation}\label{3.1}   \phi ( x) = \psi_v ( \la x , v^* \ra ) \phi (x ) \end{equation}
 for every  $ v^* \in \pi_v^k {\cal O}_v [[ t]]^{2N}$.
  If $x \notin  \pi_v ^{-k-l } t^{-1} {\cal O}_v [ t^{-1}  ]^{2N}$, we can find 
  $v^* \in \pi_v^k {\cal O}_v [[ t]]^{2N}$ such that $\psi_v ( \la x , v^* \ra )\ne 1 $, then
 (\ref{3.1}) implies  $\phi (x ) =0 $.
 \hfill $\Box$ 

\

 A function $\phi_v \in {\cal S} ( t^{-1} F_v [t^{-1} ]^{2N})$ 
   is called an elementary function if $\phi_v$ is bounded and fixed by some   
  congruence subgroup of $H_v$.  It is clear that the function $\phi_{0, v}$ in Lemma \ref{phi0padic} is an elementary function.

\begin{lemma}\label{lemma3.1}  If $v$ is a finite place, the space of elementary functions in $ {\cal S} ( t^{-1} F_v [t^{-1} ]^{2N})$
  is closed under the action of $K_v$. 
\end{lemma}

\noindent{\it Proof.}   Suppose $\phi_v$ is an elementary function, so $\phi_v$  is bounded and  fixed by 
  $ \pi_v^k {\cal O}_v (( t))^{2N} $ for some positive integer $k$.  We first note that  the $Sp_{2N} ( {\cal O}_v ((t)) )$-action
  on $H_v$ fixes $\pi_v^k {\cal O}_v (( t))^{2N} $. 
  For $ g \in K_v$ and $h\in \pi_v^k {\cal O}_v (( t))^{2N} $, we have $\bar g\in Sp_{2N} ( {\cal O}_v (( t)))$, and 
 using the identity (\ref{compatible}), we have 
 \[     h  \pi ( g ) \phi_v  =  \pi ( g )   \pi (g)^{-1}     h  \pi ( g ) \phi_v  =  \pi ( g )   ( h\cdot \bar g) v = 
    \pi ( g ) \phi_v .\]
 This proves $\pi (g) \phi_v$ is fixed by $ \pi_v^k {\cal O}_v (( t))^{2N} $.  It remains to prove
   $\pi (g) \phi_v$ is bounded.  Since $\phi_v $ and $\pi ( g) \phi_v $ are fixed by $\pi_v ^k {\cal O}_v [[ t ] ]$, 
  by Lemma \ref{lemma3.0},  $\phi_v $ and $\pi ( g ) \phi_v $ are
 supported in $ \pi_v ^{-k-l } t^{-1} {\cal O}_v [ t^{-1}  ]^{2N} $ (where $l$ is as in Lemma \ref{lemma3.0}).
  Note that $\pi ( g) \phi_v  = c T_{\bar g } \phi_v$ for some choice of Haar measure on ${\rm Im} \gamma_{\bar g}$ and some 
  scalar $c$.   Using the formula (\ref{action}) for $T_g$, we have, for $ x\in  \pi_v ^{-k-l } t^{-1} {\cal O}_v [ t^{-1}  ]^{2N} $,   

\[  | \pi ( g )  \phi ( x ) | \leq  |c|  \int_{{\rm Im} \gamma_{\bar g}} | \phi_v (  x\alpha_{\bar g} + x^* \gamma_{\bar g} ) |
      d ( x^* \gamma_{\bar g} ). \]
    Since $\phi_v $ is supported in $ \pi_v ^{-k-l } t^{-1} {\cal O}_v [ t^{-1}  ]^{2N} $ and $  x\alpha_{\bar g}
\in \pi_v ^{-k-l } t^{-1} {\cal O}_v [ t^{-1}  ]^{2N} $,
the right hand side is bounded by 
 $|c | {\rm  Vol} ( {\rm Im} \gamma_{\bar g} \cap \pi_v ^{-k-l } t^{-1} {\cal O}_v [ t^{-1}  ]^{2N} ) M $, 
 where $M$ is any upper bound of $ |\phi ( x )|$.  This proves $\pi ( g ) \phi_v $ is bounded.
\hfill $\Box$
  
\

   Let  \[ {\cal E} (  t^{-1} {\bf A} [ t^{-1}]^{2N}  ) \]
  be the space of functions on 
   $   t^{-1} {\bf A} [ t^{-1} ]^{2N} $ which are finite linear combinations of 
  $\pi ( g ) \Pi \phi_v $, where $g\in \widetilde{\Sp}_{2N} ( {\bf A} \la t \ra ) $,
   $\phi_v = \phi_{0, v}$ for all infinite places $v$ and
   almost all finite places $v$, and every remaining $\phi_v$ is an elementary function.
   It is clear that      
 ${\cal E} (  t^{-1} {\bf A} [ t^{-1}]^{2N}  )  $ is a  subrepresentation of $\otimes' {\cal S} ( t^{-1} F_v [ t^{-1} ]^{2N}  )$
  in Section 2.3. 

  We introduce a semi-subgroup of ${\rm Aut} {\bf A} \la t \ra  $ by  
\begin{eqnarray}
 & & {\rm Aut} {\bf A} \la  t \ra_{ > 1 } \nonumber \\
 & & = \{ ( \sum_{i=1}^{\infty} c_{i, v } t^i )  \in {\rm Aut} {\bf A} \la t \ra \, | \,
      |c_{1, v } | \geq 1  \, {\rm for \, all } \, v, \, \,  \,  \, \Pi_v |c_{1, v } |_v  > 1 \} \nonumber  
\end{eqnarray}
 Since ${\rm Aut} {\bf A} \la t \ra $ normalizes $\widetilde{\Sp}_{2N} ( {\bf A} \la t \ra )$, for 
    any given  ${ T} \in {\rm Aut} \bf A \la t \ra_{ > 1}   $, 
 ${T}\cdot {\cal E} ( t^{-1} {\bf A} [ t^{-1} ]^{2N} )$
  is a representation of $\widetilde {\Sp}_{2N} ( {\bf A} \la t \ra )$. 
  
We now define, for a given ${ T}\in {\rm Aut} {\bf A} \la t \ra_{ > 1 }  $,  
 a functional 
\[\theta : { T}\cdot  {\cal E} ( t^{-1} {\bf A} [ t^{-1} ]^{2N} )  \to {\Bbb C}  \]
  and prove that it is invariant under the arithmetic subgroup 
$     { Sp}_{2N} (F\la t \ra )$.
The theta functional is defined as 
 \[ \theta ( { T} \cdot f ) = \sum_{ k \in t^{-1} F[t^{-1}]^{2N} } ({ T} \cdot  f) ( k ).  \]

\begin{theorem}\label{converge} If ${ T} \in {\rm Aut} {\bf A} \la t \ra_{ > 1 }$ and $ f \in T {\cal E} ( t^{-1}{\bf A} [t^{-1} ]^{2N} )$,
 then 
 $\theta ( f )$ converges absolutely .
 \end{theorem}

  This theorem is the adelic version of the Theorem 4.9 in \cite{Z}, we sketch its proof.
First we need some lemmas analogous to Lemma 4.7 and 4.8 in \cite{Z}.

\begin{lemma}\label{bipos}  Let $S$ be a finite  set of places of $F$ that contains all the infinite places.
  If $f\in {\cal S } ( \Pi_{v\in S} F_v )$ is bipositive (i.e. $f$ and its Fourier transform 
   $Ff$ satisfy $f\geq 0$ and $Ff \geq 0 $).  Suppose $L\subset  \Pi_{v\in S} F_v$
  is a lattice (i.e. $L$ is a discrete subgroup with compact quotient of $ \Pi_{v\in S} F_v$),
  then for all $v\in  \Pi_{v\in S } F_v$,
 \[   \sum_{n\in L}  f( n + v  ) \leq    \sum_{n\in L}  f( n  ). \]
\end{lemma}

\noindent{Proof.}   Apply the Poisson summation formula, we have
  \[  \sum_{n\in L}  f( n + v  ) = \sum_{n\in L'} \psi(  (v , n ) ) Ff ( n ) , \]
where $L'$ denotes the dual lattice of $L$.
  So 
\[  \sum_{n\in L}  f( n + v  )\leq \sum_{n\in L'} | \psi(  (v , n ) ) Ff ( n ) |
   = \sum_{n\in L'}  Ff ( n ) = \sum_{n\in L}  f ( n ).  \]
\hfill $\Box$

\

  A multi-dimensional generalization of the lemma is

\begin{lemma}\label{mbipos} 
  Let $S$ be as in Lemma \ref{bipos}.  If $f_i \in {\cal S} ( \Pi_{v\in S } F_v )$ ($i=1 , \dots , m$
  are bi-positive functions,
 and $f = \Pi_{i=1}^m  f_i$ be regarded as a function on $ {\cal S} ( \Pi_{v\in S } F_v^m  )$
 in the obvious way.  If  
  $ N_v $ for each $v\in S $ is an $m\times m$ upper triangular unipotent matrix acting on $F_v^m$,
  let $ N = \Pi_v N_v $ be the product of the $N_v$ acting on $\Pi_{v\in S } F_v^m $ and  
  let ${\cal O}_S$ be the ring of $S$-integers that is embedded in $\Pi_{v\in S } F_v$ diagonally,
so $L = {\cal O}_S^m  $ is a lattice in $\Pi_{v\in S } F_v^m $.  Then for any $c\in F$,
\[ \sum_{n\in c L}  f( N n   ) \leq    \sum_{n\in c L}  f( n  ). \]
\end{lemma}

\noindent The proof of this lemma is similar to that of Lemma 4.8 in \cite{Z}, it uses Lemma \ref{bipos}.

\

\noindent{\it Proof of Theorem \ref{converge}.}      For simplicity, we assume
   $ { T} =\{{ T}_v \} $ with ${T}_v = 1 $ for all finite places $v$.  
  For each infinite place $v$, ${ T}_v$ can be factorized as a product 
 $ { T}_v = {T}_{v, d } {T}_{v, u} = ( q_v t ) \circ ( t+  \sum_{i=2} c_{i } t^i )$, and we assume $|q_v | > 1 $. We may 
 assume $ f = \Pi f_v$ such that $f_v =\phi_{v, 0}$ for almost all finite places $v$. 
  Let $ S$ be a finite set of places containing all infinite $v$ and all finite $v$ with 
  $f_v \ne \phi_{v , 0 }$. Then
 \[  \theta ( f ) =  \sum_{k \in t^{-1} {\cal O}_S [t^{-1}]^{2N} } ({ T} \cdot  \Pi_{v\in S} f_v ) ( k ). \]    
  For each $v\in S$, $ f_v = g_v \phi_v $ for $g_v \in \widetilde{Sp}_{2N} ( F_v (( t)))$
  and $\phi_v$  elementary for a finite place $v$ and $ \phi_v = \phi_{v, 0 }$ for an infinite place $v$.
  By the Iwasawa decomposition (\ref{I1}), we write 
 $g_v =  b_v   k_v $ for $k_v \in K_v$, $b_v \in B_v$, and we further
  write $b_v= u_v a_v$, where $u_v$ is in the pro-unipotent radical of $B_v$ and $a_v$ is in the Cartan subgroup.
   By Lemma \ref{lemma3.1}, for a finite place $v$,
  $k_v f_v$ is again elementary, and it is easy to see that $a_v k_v f_v $ is elementary.
 So we may assume $f_v = u_v \phi_v $ for some 
   elementary function $\phi_v$.   Since $\phi_v $ is elementary,  
   $\phi_v$ is bounded by  a constant multiple $C_v$ of the characteristic function $\phi_v'$ of 
  $ \pi_v^{k_v} t^{-1} {\cal O}_v [t^{-1}]^{2N}$ for some integer $k_v$. 
  Then  
\[ |u_v \phi_v ( x ) | <   C_v | u_v \phi_v' ( x ) |. \]    
So we may further assume that  $\phi_v$ is the characteristic function of $ \pi_v^{k_v} t^{-1} {\cal O}_v [t^{-1}]^{2N}$.
For an infinite place $v$, 
  $f_v =  { T}_v g_v   \phi_{v, 0 } $.   We write $g_v = u_v a_v k_v $, so 
    \[ f_v = { T}_{v, d } { T}_{v, u } u _v a_v \phi_{v, 0 }   =  { T}'_{v, u } { T}_{v, d } u _v a_v \phi_{v, 0 },  \]
  where $ { T}'_{v, u } =  { T}_{v, d } { T}'_{v, u }{ T}_{v, d }^{-1}$.
So we need to prove the convergence of 
\begin{equation}\label{t1}
   \sum_{k \in t^{-1} {\cal O}_S [ t^{-1} ]^{2N}   }  
 \Pi_{v\in S_{\infty}} \, | \, { T}'_{v, u } \, ((q_v t)\cdot u_v) \,  (q_v t )\,  a_v \phi_{v,0} ) ( k )|
   \Pi_{v\in S_{\rm fin} } | u_v \phi_v ( k ) | .  \end{equation}
We first consider the sum over the finite dimensional subset  $  {\cal O}_S [ t^{-1} ]^{2N}_d t^{-1}$ of  
$  {\cal O}_S [ t^{-1} ]^{2N} t^{-1}$  consisting
  of polynomials of degree less or equal to $d$.   Since $  (q_v t )\,  a_v\phi_{v, 0}$ ($v\in S_{\infty}$) is a Gaussian 
 function, its restriction on every finite dimensional subspace of $ t^{-1} F_v [ t^{-1} ]^{2N} $ 
 is bipositive; and it is also easy to see that the restriction of $  \phi_v $ ($v\in S_{\rm fin}$) 
 on every finite dimensional subspace of $ t^{-1} F_v [ t^{-1} ]^{2N} $ is bipositive. And 
   the operators  ${T}_{v, u }  ((q_vt)\cdot u_v)$ ($v\in S_{\infty}$) and
  and $u_v$ ($v \in S_{\rm fin}$) have unipotent matices under the standard monomial basis of $  F_v [ t^{-1} ]^{2N}_d t^{-1}$.
  Therefore we may apply Lemma \ref{mbipos} to prove that
 \[ \sum_{k \in t^{-1} {\cal O}_S [ t^{-1} ]_d^{2N}   }  
 \Pi_{v\in S_{\infty}} \, | \, { T}_{v, u } \, ((q_v t)\cdot u_v) \,  (q_v t )\,  a_v \phi_{v,0} ) ( k )|
   \Pi_{v\in S_{\rm fin} } | u_v \phi_v ( k ) | \]
 is bounded  above by 
  \[   \sum_{k \in t^{-1} {\cal O}_S [ t^{-1} ]_d^{2N}   }  
 \Pi_{v\in S_{\infty}} \, |    (q_v t )\,  a_v \phi_{v,0} ) ( k )|
   \Pi_{v\in S_{\rm fin} } |  \phi_v ( k ) | . \]
Letting $d \to \infty $, we see that  
 \[     {\rm (\ref{t1} )} \leq  
      \sum_{ k \in  \frac 1 M t^{-1} {\cal O}  [t^{-1} ]^{2N}}
     \Pi_{v\in S_{\infty}} ( (q_v t) a_v \phi_{v, 0} )( k) , \]
   where $M$ is an integer (with prime divisors only  in $S$, above) depending only on 
 $\otimes_{v\in S_{\rm fin} } \phi_v $ and ${\cal O} $ is the ring of integers in $F$ .  
  Using 
 the fact that  $ q_v > 1 $ for all $v\in S_{\infty}$, we can prove the right hand side is convergent.
  \hfill $\Box $

 \

We expect that if  $ T = ( \sum_{i=1}^{\infty} c_{i, v} t^i ) \in {\rm Aut} {\bf A} \la t \ra $ in Theorem \ref{converge} 
 satisfying the weaker condition  $\Pi_v | c_{1, v} | > 1 $, then Theorem \ref{converge} still holds.

\begin{theorem}\label{inv}  If 
   $g\in {\rm Sp}_{2N} ( F \la  t \ra )$,  then $\theta ( g f  ) = \theta ( f ) $ for $f$ as in Theorem \ref{converge}. 
 \end{theorem}

 \noindent  {\it Proof.}  Recall that $ {\rm Sp}_{2N} ( F\la  t \ra )$ is generated by root vectors 
  $x_{\alpha } ( a ) $ for $\alpha \in \Delta_0$ and $a \in F\la t \ra $. 
   If $a \in F [[ t]]$, it is clear that $\theta $ is invariant under $ x_{\alpha } ( a )$.
  The general $ \pi ( x_{\alpha } ( a )) $ can be written as 
   $ T_{t^m}   \pi ( x_{\alpha  } ( b ) )  T_{t^m}^{-1} $ (see (\ref{2a1})) for some $b \in F[[t]]$.   Since $ T_{t^m} $ acts as 
 partial Fourier transform, it preserves $\theta $ by Poisson summation formula. Therefore $  \pi ( x_{\alpha } ( a )) $ preserves $\theta $.
\hfill $\Box $

\

By Theorem \ref{act}, we know that for $g \in \widetilde{\Sp}_{2N} ( {\bf A} \la t \ra  ) $,
and  $f\in { T}\cdot  {\cal E} ( t^{-1} {\bf A} [ t^{-1} ] ^{2N} )$, $ \pi ( g ) f $ is a scalar multiple of  $T_g f$, 
  where $T_g$ is given by (\ref{action}).   Using the theta functional, we can determine this scalar
  in the case $g\in {\Sp}_{2N} ( F\la t \ra )$.

\begin{prop}\label{prop}  For $g\in {\Sp}_{2N} ( F\la  t \ra )$, $f$ as above, then $\pi(g)$ is equal to 
 $T_g$ given by 
    \[ ( T_g f ) ( x ) 
= \int_{{\rm Im} \gamma_g ({\bf A})}    S_{ g}( x + x^* )   f ( x \alpha + x^* \gamma ) d ( x^*\gamma  ) \]
where 
\[ S_{g} ( x + x^*) = \psi \left( \frac 12 \la x\alpha , x \beta \ra +
          \frac 12 \la x^* \gamma , x^* \delta \ra + \la x^* \gamma , x \beta \ra \right) ;\]
  with the Haar measure given by the 
  condition that the covolume of  ${\rm Im} \gamma_{g} (F) $
  in ${\rm Im} \gamma_g ({\bf A})$ 
 is $1$
 (In particular if ${\rm Im } \gamma_g  =\{ 0 \} $, ${\rm Im } \gamma_g (\bf A )  $ is a point, its volume is taken as $1$) . 
\end{prop}

\noindent {\it Proof.}  By Theorem \ref{inv},  $\pi ( g ) $ preserves 
  the theta functional, and we also known that $ \pi ( g) = c_g T_g $ for some scalar $c_g$. 
  To prove $\pi ( g ) = T_g$, it is sufficient to prove $T_g$ also 
  preserves the theta functional. 
  To prove this, we use the Bruhat decomposition of 
 \[ 
\Sp_{2N} ( F \la t \ra ) = \sqcup_w  \Sp_{2N} ( F \la t \ra_+  ) w \Sp_{2N} (  F \la t \ra_+  ) ,\]
where $ F\la t \ra_+ = F\la t \ra \cap {\bf A} \la t \ra_+ $,
 and $w$ runs through all diagonal matrices $t^m $.   It is easy to verify that 
  $T_g $ for $g\in  Sp_{2N} ( F\la t \ra_+ )$ preserves the theta functional, 
  and $T_{t^m}$ is a partial Fourier transform, so it also preserves the theta functional.
   To prove $T_g$ preserves the theta functional, we write
  $ g = g_1 w g_2$ according to the Bruhat decomposition.  Using Lemma 2.5 \cite{Z},
 we have $T_g = T_{g_1} T_w T_{g_2} $. This proves the proposition.
 
\

\noindent \section{  Eisenstein Series for Loop  Metaplectic Groups  }

\noindent   We assume $F$ is a number field.   Let $ F^{2n}$ be  the standard symplectic
  space and $V$ be an $m$-dimensional $F$-space with a non-degenerate symmetric,  ani-isotropic, bilinear form 
  $(\, , \,): V\times V \to F$ . Then  $F^{2n}\otimes V $
  is an
   $F$-symplectic space with the symplectic form 
\[ \la v_1\otimes m_1 , v_2 \otimes m_2 \ra = \la v_1 , v_2 \ra ( m_1 , m_2 ) . \]   
   Let $e_1 , \dots , e_{2n}$ be the standard basis for $F^{2n}$, we have 
   $ \la e_i , e_{i+n} \ra = 1 $ for $i=1, \dots , n$ and all other symplectic pairings are $0$.
  Let $f_i$ be a basis of $V$, and $f_i'$ be the dual basis of $V$, i.e. 
  $( f_i , f_j') = \delta_{i, j}$.  Then the $e_i \otimes f_j $ ($ 1\leq   i \leq n , 1\leq  j  \leq m$)   
  and the $ e_{n+i} \otimes f_j'$ ($ 1\leq   i \leq n , 1\leq  j  \leq m$) together form a symplectic basis of $F^{2n} \otimes V$.
 We fix this choice of symplectic basis, and  therefore identify the symplectic space $F^{2n} \otimes V$
  with the standard one $ F^{2N}$ ($N = m n $).  
  We shall define  Eisenstein series for $\widetilde {\rm Sp}_{2n} ( {\bf A} \la t \ra ) $
  for a function in $ \otimes'{\cal S} ( t^{-1} F_v [ t^{-1} ]^{2N} )$   
 satisfying certain conditions.  We then prove the convergence of these Eisenstein series reduces to the convergence of 
 Eisenstein series for a certain height function on $\widetilde {\rm Sp}_{2n} ( {\bf A} \la t \ra ) $ (Lemma 4.4).

    We have an embedding
  \[ {\rm Sp}_{2n} \times G \subset {\rm Sp}_{2N}  ,\]
 where $G$ denotes the orthogonal group of $V$,  and an embedding 
   \[ {\rm Sp}_{2n} ( F(( t))) \times  G ( F((t))) \subset {\rm Sp}_{2N} ( F((t))). \]

 For a place $v$ of $F$,  we have 
  the Weil representation  
  of $  \widetilde{Sp}_{2N} ( F_v (( t )))$ on  $ {\cal S} (  F_v^{2n}  \otimes V_v [t^{-1} ] t^{-1} )$   as in Section 2.1.
  The embedding $ {\rm Sp}_{2n} \to {\rm Sp}_{2N}$ gives a homomorphism 
  $ \widetilde {\rm Sp}_{2n} ( F_v (( t)) )  \to \widetilde {\rm Sp}_{2N} ( F_v (( t)))$,
  under which an element $c$ in the center ${\Bbb C}^* \subset \widetilde {\rm Sp}_{2n} ( F_v (( t)) ) $,
   goes to $ c^m $ in the center  ${\Bbb C}^*$  of $ \widetilde {\rm Sp}_{2N} ( F_v (( t)) ) $.  For almost all $v$, 
 the ``maximal compact" subgroup $K_v$ of $\widetilde {\rm Sp}_{2n} ( F_v (( t)) )$   
   maps into the ``maximal compact" subgroup of $\widetilde {\rm Sp}_{2N} ( F_v (( t)) )$.
  Therefore we have a morphism $  \widetilde{\rm Sp}_{2n} ( {\bf A} \la t \ra ) 
  \to  \widetilde{\rm Sp}_{2N} ( {\bf A} \la t \ra ) $. By Lemma \ref{Lemma2.4a},  
 we may view ${ G }( F_v [[t ]]) \subset { G }( F_v ((t )))$ as a subgroup of $\widetilde{ \rm Sp }_{2N} ( F_v ((t)))$,
 so we may view  ${ G }( {\bf A}  \la t \ra_+)$ 
  as a subgroup of $\widetilde{ \rm Sp }_{2N} ( {\bf A} \la t \ra )$.

\begin{lemma}\label{Lemma4.1}  The actions of $\widetilde{Sp}_{2n} ( F_v ((t)))$ and $ G ( F_v [[t ]] )$  commute  .
\end{lemma}

 This result is proved directly using (\ref{action}).

\

For each $\phi = \Pi \phi_v  \in {\cal E} ( t^{-1} {\bf A} [ t^{-1}]^{2N} )  $,  and $ T \in {\rm Aut} {\bf A} \la t \ra_{> 1 } $,
 we consider  the  function $\Phi $ on $  T   \widetilde {\rm Sp}_{2n}( {\bf A} \la t \ra )$ by 
  \begin{equation}\label{4.1}   \Phi ( T  g )   =    ( T  g \phi ) ( 0 ) . \end{equation}
Note that  $ \widetilde {\rm Sp}_{2n} ( {\bf A} \la t \ra  )$ is normalized by $T$,  
the set $T  \widetilde{\rm Sp}_{2n} ( {\bf A} \la t \ra )$ is stable under 
  both left and right multiplication by $ \widetilde {\rm Sp}_{2n} ( {\bf A} \la t \ra  )$. And it is clear that 
  $\Phi (  h T g ) = \Phi ( T g )$ for 
   $ h\in  {\rm Sp}_{2n} ( {\bf A} \la t \ra_+ )$. We define 
the Eisenstein
  series 
 \begin{equation}\label{Eisen} {\rm E} (\phi ,    T g ) =   \sum_{ r \in  {\rm Sp}_{2n} ( F\la t \ra_+ ) \setminus {\rm Sp}_{2n} ( F\la t \ra )} 
   \Phi ( r T g ) . \end{equation}

 Although our results can be proved for more general $T$, 
for simplicity of exposition, we make the following assumption on $T$: 

\noindent {\bf Assumption 4.1.}  $T = ( T_v ) \in {\rm Aut } {\bf A} \la t \ra $ satisfies the following conditions:
   each $p$-adic component $T_v$  is the identity element of the group $ {\rm Aut}  F_v (( t))$,
  and each infinite component $T_v$  is $ q_v t $ with $q_v \in F_v $ such that $  |q_v| > 1 $.  

  The main result of this section is a comparison lemma (Lemma 4.4) about the Eisenstein series (\ref{Eisen}) with 
  the  Eisenstein series defined by a certain height function. 
 For $T$ as in Assumption 4.1,
    we define for each place $v$, the local height function
  $ \tilde h_v :  T_v \widetilde{\rm Sp}_{2n} ( F_v (( t))) \to {\Bbb C}^* $ as follows,  for $g_v\in \widetilde{\rm Sp}_{2n} ( F_v (( t)))$,
  we write $g_v = c b k $ for $ k \in K_v$, $b\in Sp_{2n} ( F_v [[ t]] )$, $c \in {\Bbb C}^*$, we define $\tilde h_v ( T_v g_v )= | c |$.  We define 
  for $g= \Pi g_v  \in  \widetilde{\rm Sp}_{2n}( {\bf A} \la t \ra )$,   $\tilde h ( T g ) = \Pi_v \tilde h_v ( T_v g_v )$.

\

 Both $\Phi $ and $ \tilde h^m $ are left invariant under
$ Sp_{2n} ( {\bf A} \la t \ra_+ )$, and they are equal on the center $\oplus_v {\Bbb C}^*$. 
To compare them, it is enough to compare their restrictions on $ \Pi K_v $.  
In the classical case, since $\Pi_v K_v$ is compact, we have  
$  \Phi ( g ) \leq C \tilde h (g)^m $ 
for some scalar $C$, so  the convergence
  of the Eisenstein series associated to $\Phi $ reduces to the Eisenstein series associated to $\tilde h^m$, where the Godement criterion can be applied.
  In our case, the argument that uses 
  the compactness of $ \Pi K_v$ doesn't apply, we need to assume more conditions on $\phi $ to prove the comparison inequality involving $\Phi$ and $\tilde h^m $.

   \

To give our conditions on $\phi $, we use  the action of the local Heisenberg group  $ F_v (( t))^{2N} \times F_v$  on 
  $ {\cal S} ( t^{-1} F_v [ t^{-1} ]^{2N} )$ (Section 2.1).   Assume the local additive character $\psi_v$ that is used
  to define the representation of the Heisenberg group has conductor ${\pi_v}^{n_v} {\cal O}_v $ (i.e. $\psi_v ( \pi_v^{n_v} {\cal O}_v ) = 1$
 and $\psi_v ( \pi_v^{n_v-1 } {\cal O}_v ) \ne  1$).  
   If $2k\geq n_v $,  then 
   action of any two $ a , b \in  \pi_v^{k} {\cal O}_v (( t))^{2N} $ commutes.   
  If $\phi_v \in   {\cal S} ( t^{-1} F_v [ t^{-1} ]^{2N} )$ is fixed by   $ \pi_v^{k} {\cal O}_v (( t))^{2N}$ for 
  some $k \geq  n_v / 2  $, then in particular $ \phi_v $ is fixed by $  \pi_v^{k} t^{-1}{\cal O}_v [t^{-1}]^{2N}$, so 
  $\phi_v$ is constant on each coset of  $ t^{-1} F_v [ t^{-1} ]^{2N} $ mod $ \pi_v^{k} t^{-1}{\cal O}_v [t^{-1}]^{2N}$.
  Since $\phi_v $ is also fixed by $  \pi_v^{k} {\cal O}_v [[t]]^{2N}$, then $\phi_v $ is supported on
   $ t^{-1} \pi_v^{n_v -k } {\cal O}_v [ t^{-1} ]^{2N} $.  Therefore we may regard $\phi_v $ as function on 
  \[ t^{-1}  ( \pi_v^{n_v - k } {\cal O}_v / \pi_v^k {\cal O}_v ) [ t^{-1} ]^{2N} .\]  
   For any non-negative integer $k$ with $2k \geq n_v$, we let 
  
\[  {\cal S}_{v,k} = \{  f \in {\cal S}  ( t^{-1} F_v [ t^{-1} ]^{2N} ) \, | \, f \, {\rm is \, fixed \, by \, } \, \pi_v^{k} {\cal O}_v (( t))^{2N} \}.  \] 
  We denote the ``maximal compact" subgroup of $ \widetilde{\rm Sp}_{2N} ( F_v (( t)))$ by $K_{v, N}$.

\begin{lemma}  ${\cal S}_{v,k}$ is closed under the action of $K_{v, N}$, and therefore closed under the action of $K_v$.
\end{lemma}

 This lemma follows form the proof of Lemma 3.2.  The proof uses the fact that $K_{v, N}$ normalizes  $\pi_v^{k} {\cal O}_v (( t))^{2N} $.
 As we remarked earlier, a function  $f\in {\cal S}_{v,k} $ can be viewed as an element in 
 {\ \[{\bf C} \left( t^{-1}  ( \pi_v^{n_v - k } {\cal O}_v / \pi_v^k {\cal O}_v ) [ t^{-1} ]^{2N} \right), \]
  the space of complex valued  functions on $t^{-1}  \left( \pi_v^{n_v - k } {\cal O}_v / \pi_v^k {\cal O}_v \right) [ t^{-1} ]^{2N} $.
 Let  $f \mapsto \hat f $ denote the isomorphism from $ {\cal S}_{v,k}$ to $ {\bf C} ( t^{-1}  ( \pi_v^{n_v - k } {\cal O}_v / \pi_v^k {\cal O}_v  ) [ t^{-1} ]^{2N} )$.
 We define 
 \[ {\cal S}_{v,k, c}= \{ f \in {\cal S}_{v, k } \, | \, \hat f \, {\rm \, has \, finite \, support .}  \} \]
 It is clear that  the space  $ {\cal S}_{v,k, c}$ is isomorphic to 
 \begin{equation}\label{c0}  {\bf C}_0 \left( t^{-1}  ( \pi_v^{n_v - k } {\cal O}_v / \pi_v^k {\cal O}_v ) [ t^{-1} ]^{2N} \right), \end{equation}
   the space of complex valued  functions on $t^{-1}  ( \pi_v^{n_v - k } {\cal O}_v / \pi_v^k {\cal O}_v ) [ t^{-1} ]^{2N} $ with finite support.
 
 We introduce an inner product on ${\cal S}_{v,k, c}$. For $f_1 , f_2 \in {\cal S}_{v,k, c}$, 

\begin{equation}\label{inner}  (f_1 , f_2 ) = \sum_{x} {\hat f}_1 ( x )\overline {{\hat f}_2  ( x )}, \end{equation} 
where the sum is over $x\in t^{-1}  ( \pi_v^{n_v - k } {\cal O}_v / \pi_v^k {\cal O}_v ) [ t^{-1} ]^{2N}$, since $\bar f_1$ and $\bar f_2$ have
  finite support, (\ref{inner}) is a finite sum.

\
 
\begin{lemma}  ${\cal S}_{v,k, c}$ is closed under the action of $K_{v, N}$ and $K_v$, and the actions of 
  $K_{v, N}$ and $K_v$ are unitary.
\end{lemma}

Lets give a conceptual explanation of this lemma.  We first use the additive character $\psi_v $ to define a bi-character of 
 the finite abelian group $ \pi_v^{n_v - k } {\cal O}_v / \pi_v^k {\cal O}_v $ as follows.  
  For any its elements $a, b$, let $\bar a , \bar b \in \pi_v^{n_v - k } {\cal O}_v$ be liftings.
  Since $\psi_v$ has conductor $\pi_v^{n_v} {\cal O}_v$,  $\psi_v ( \bar a \bar b)$ depends only on $a, b$, not on the choice of their liftings.
  We denote this bi-character by $\psi ( a b)$.  Using this bi-character and the symplectic structure on $\pi_v^{n_v} {\cal O}_v (( t))^{2N}$,
  we have a non-degenerate skew symmetric bi-character which we denoted by $\psi ( f , g )$ on the locally compact abelian group 
   \[   (\pi_v^{n_v - k } {\cal O}_v / \pi_v^k {\cal O}_v ) (( t))^{2N}.  \]
(the topology is the $t$-adic topology).  The subgroups $ (\pi_v^{n_v - k } {\cal O}_v / \pi_v^k {\cal O}_v ) [[ t ]]^{2N}$ (which is compact)
  and $ t^{-1} (\pi_v^{n_v - k } {\cal O}_v / \pi_v^k {\cal O}_v )  [ t^{-1} ]^{2N}$ (which is discrete) are complementary maximal isotropic subgroups.
 We have the associated 
  Heisenberg group 
 \begin{equation}\label{reduced}
  \bar H = (\pi_v^{n_v - k } {\cal O}_v / \pi_v^k {\cal O}_v ) [[ t ]]^{2N}\times  t^{-1} (\pi_v^{n_v - k } {\cal O}_v / \pi_v^k {\cal O}_v )  [ t^{-1} ]^{2N} \times S^1
 \end{equation}
  with the group law given by the relations
\[  f g = \psi ( f , g ) g f  \]
  for $f\in \pi_v^{n_v - k } {\cal O}_v / \pi_v^k {\cal O}_v ) [[ t ]]^{2N}$ and $g \in  t^{-1} (\pi_v^{n_v - k } {\cal O}_v / \pi_v^k {\cal O}_v )  [ t^{-1} ]^{2N} $.
  Clearly $Sp_{2N} ( {\cal O}_v (( t)))$ acts on $\bar H$ as automorphisms, and the action  factors through
   $Sp_{2N} ( ({\cal O}_v/ \pi_v^{2k- n_v} {\cal O}_v ) (( t)))$. 
  Note also that the group $\bar H$ is a locally compact topological group. By the Stone-Von Neumann Theorem,
 there is a unique (up to isomorphism) irreducible unitary representation of $\bar H$ with the central character $S^1 \to S^1 , z \mapsto z$.
  The space $ L^2 (  t^{-1} (\pi_v^{n_v - k } {\cal O}_v / \pi_v^k {\cal O}_v )  [ t^{-1} ]^{2N} )$ is a model of this representation, where 
  the elements in $t^{-1} (\pi_v^{n_v - k } {\cal O}_v / \pi_v^k {\cal O}_v )  [ t^{-1} ]^{2N}$ acts as translations and the elements in 
$(\pi_v^{n_v - k } {\cal O}_v / \pi_v^k {\cal O}_v ) [[ t ]]^{2N}$ acts as multiplication by additive characters. 
  The smooth vectors of this representation is precisely (\ref{c0}) with the inner product given by (\ref{inner}).
  By the uniqueness, certain central extension of 
 $Sp_{2N} ( {\cal O}_v (( t)))$ acts on  $L^2 (  t^{-1} (\pi_v^{n_v - k } {\cal O}_v / \pi_v^k {\cal O}_v )  [ t^{-1} ]^{2N} )$ and therefore acts on (\ref{c0}).
  This representation is the one in Lemma 4.3. 

 We can also prove Lemma 4.3 directly.
  We use the fact that $K_{v, N}$ is generated by ${\rm Sp}_{2N} ( {\cal O}_v [[t]] )$ 
and $K_{\alpha_0}$, which is the maximal compact 
 subgroup of the metaplectic group for $SL_2 ( F_v )$ associated to the extra simple root 
in the affine root system (see Section 2.2).  It is easy to verify that $Sp_{2N}  ( {\cal O}_v )$ preserves 
  $ {\cal S}_{v,k, c}$ and is unitary.  By a direct computation, we prove that $K_{\alpha_0}$ preserves   $ {\cal S}_{v,k, c}$ and is unitary. Similar proofs works
  for $K_v$.

\begin{lemma}\label{compare2} If $\phi =\Pi_v \phi_v$ satisfies the conditions that (1) for each infinite place $v$, $\phi_v$ is $\phi_{v, 0}$ as in Lemma 2.5 
  (2) for all but finitely many finite places $v$, $\phi_v $ is $\phi_{v, 0}$ as in Lemma 2.6, (3) each of the remaining components $\phi_v $
   is in ${\cal S}_{v,k, c}$ for some $k$.  Let $\Phi $ be as (\ref{4.1}). Then 
  there is a constant $C$ depending on $\phi$ only, such that  
 \[    | \Phi ( T g  ) | \leq C \cdot \tilde h  ( T g  )^m    \] for all $g\in    \widetilde{\rm Sp}_{2n} ( {\bf A} \la t \ra  )$.
\end{lemma}

   \noindent {\it Proof.}  Since $T$ normalizes ${\rm Sp}_{2n} ( {\bf A} \la t \ra_+)$,  and $\Phi ( b T g )= \Phi ( T g ) $,
  $\tilde h ( b T g )^m = \tilde h ( T g )^m $ for $b \in \widetilde {\rm Sp}_{2n} ( {\bf A} \la t \ra_+ )$. And for  $c$ in the center
  $ {\Bbb C}^*$,  $\Phi ( c T g )=  | c|^{m} \Phi ( T g )$, $\tilde h ( c T g )^m = |c|^m  \tilde h  ( T g )^m $.
 It is enough to prove the inequality for $g \in \Pi_v K_v $.
  By the conditions on $\phi$, it is suffices to prove that for each exceptional place $v$ in (3), 
   the function $ K_v \to {\Bbb C} $ given by $ k \mapsto ( k \cdot \phi_v ) ( 0 ) $ is bounded. It is proved using Lemma 4.3:
\[  | ( k \cdot \phi_v )  ( 0 ) |^2  \leq  ( \widehat { k \cdot \phi_v } , \widehat { k \cdot \phi_v } ) =  (\hat \phi_v , \hat \phi_v ). \]
    \hfill $\Box$

\

  Lemma \ref{compare2} implies that the Eisenstein series (\ref{Eisen}) is majorized by $C$ times the Eisenstein series
\begin{equation}\label{Eisen2} 
 \tilde E_m ( T g ) =   \sum_{ r \in  Sp_{2n} ( F[[t ]] ) \setminus Sp_{2n} ( F(( t )))} 
    \tilde h  ( r T g )^m  , \end{equation}  
  the convergence of which is proved  in the next section.

%
%
%
%
%

\section{Convergence of Eisenstein Series}

In this section, we prove the convergence of Eisenstein series associated to 
 the ``Siegel parabolic" subgroup of the symplectic loop group defined using the inverse tame symbol. The convergence
 result implies the convergence of (\ref{Eisen2}). We follow the method used in \cite{Garland1}  \cite{Garland2}, where the convergence is proved for
 Eisenstein series  induced from a Borel subgroup.  

Let $F$ be a number field and $v$ be a place of $F$. Since the inverse tame symbol on $F_v (( t ))$ satisfies 
 the Steinberg relations, it gives a central extension of $\Sp_{2n} ( F_v(( t )) ) $, which we 
 denote by  $\widehat {\Sp}_{2n} ( F_v (( t ))) $:

\[  1\to  F_v^* \to \widehat {\Sp}_{2n} ( F_v (( t )) ) \to {\Sp}_{2n}( F_v  (( t )) ) \to 1 .\]
The group  $\widehat {\Sp}_{2n} ( F_v (( t )) )$ is generated by root vectors $x_{\alpha } ( a ) $ with relations
 (\ref{llg1}), (\ref{llg2}), and (\ref{llg5}) with the symbol replaced by $ C( a , b)^{-1}$.
 We can similarly define Borel subgroup $B_v$ and ``maximal compact" subgroups  $K_v$ of  $\widehat {\Sp}_{2n} ( F_v (( t ))) $ as in Section 2.2.
 The Iwasawa decomposition holds:
 \[  \widehat {\Sp}_{2n} ( F_v (( t )) ) = B_v K_v , \]
 see \cite{S}. 
 Let  ${\rm Aut } \, F_v ((t ))$ be the  reprametrization group given in Section 2.2.  It acts on
  $ {\Sp}_{2n}( F_v  (( t )) )$  by changing  the variable $t $, i.e., for $\sigma ( t ) \in {\rm Aut } \, F_v ((t ))$, \
  $ g ( t)  \in  {\Sp}_{2n}( F_v  (( t )) )$, 
  \begin{equation}\label{5.1}
  \sigma ( t ) \cdot g ( t)  = g ( \sigma ( t ) ) . \end{equation}
 Since the changing variable action preserves the tame symbol, it has a compatible  action on 
  $\widehat {\Sp}_{2n} ( F_v (( t )) )$ given by 
   \[  \sigma ( t ) \cdot x_{\alpha } ( a (t ) )  = x_{\alpha } ( a  ( \sigma  ( t )   )) , \]
 and the center $F_v^*$ is fixed. 
  We shall only consider the subgroup 

 \[   \sigma ( F_v^* t )   = \{ q t   \, | \, q \in F_v^* \} \subset {\rm Aut} \, F_v (( t)) .   \]
 We have the semi-direct product group 

 \[   \widehat {\Sp}_{2n} ( F_v (( t )) )  \rtimes  \sigma ( F_v^* t  ) . \]

 \noindent Let ${\rm T}_v$ denote the subgroup   generated by $ h_{\alpha } (c ) $ ($\alpha \in \Delta_0$, $c \in F_v^*$),
  and the center  $F_v^*$.
Then ${\rm T}_v  $ is a ``maximal torus" of $\widehat {\Sp}_{2n} ( F_v (( t )) )$, which is isomorphic to $(F_v^*)^{n+1}$ . 
 Let 
$ \Pi_v \widehat {\Sp}_{2n} ( F_v (( t )) ) $ be the direct product.  It has the center $ \Pi_v F_v^* $ which contains the idele group 
 ${\bf A}^*$.
 We denote 
  \[  \widehat {\Sp}_{2n} ( {\bf A} \la  t \ra  ) \deff  \Pi_v' \widehat {\Sp}_{2n} ( F_v (( t )) ) \]
 where $\Pi_v'$ is 
the restricted product relative the $K_v$'s.  Note that 
 we have the exact sequence 
 \[ 1\to  {\bf A}^* \to    \widehat {\Sp}_{2n} ( {\bf A}  \la  t \ra )           \to  {\Sp}_{2n} ( {\bf A} \la  t \ra ) \to 1 ,\] 
    In particular $F^* \subset {\bf A}^* $ is in center of 
 $ \widehat {\Sp}_{2n} ( {\bf A} \la  t \ra )$. 
 \noindent We have an embedding $\Sp_{2n} ( F\la  t \ra) \to   \widehat {\Sp}_{2n} ( {\bf A}\la t \ra   )  / F^*  $ given by 
 \[ x_{\alpha } ( a) \mapsto \Pi_v x_{\alpha } ( a ) .\]

\noindent The group $ \sigma  ( {\bf A}^* t  ) = \Pi_v' \sigma  (F_v^* t )$ acts on 
 $ \widehat {\rm Sp}_{2N} ( {\bf A } \la t \ra  )$, the action being  induced from the local actions.  We have semi-direct product group
 $   \widehat  {\rm Sp}_{2N} ( {\bf A } \la t \ra  ) \rtimes \sigma  ( {\bf A}^* t  ) $,  and we set
 ${\rm T}_{\bf A} = \Pi_v' {\rm T}_v \cong {{\bf A}^*}^{n+1} $.

 We recall the affine Kac-Moody algebra $\hat {\frak g} $ for a  complex simple Lie algebra $\frak g$:

 \[  \hat {\frak  g} = \frak g \otimes {\Bbb C} [ t , t^{-1} ] + {\Bbb C } K + {\Bbb C} d .\] 

\noindent The Lie bracket is given by 

 \[  [ a \otimes t^m  ) , b \otimes  t^n  ) ]  = [ a , b ] \otimes t^{m+n} +  ( a , b ) \delta_{m+n, 0 }   K  ,\] 
  \[   [ d,  a \otimes t^n ] = n a \otimes t^n  ,\]
and $K$ is in the center. The bilinear form $(a , b)$ above is the normalized Killing form such that $( \check \theta ,   \check \theta ) = 2$,
 where $ \check \theta $ is the coroot corresponding to the longest root $\theta $. 
Let $ \frak h$ be a Cartan subalgebra of $\frak g$, then 
\[ \hat {\frak h} =\frak h + {\Bbb C} K + {\Bbb C} d \]
 is a Cartan subalgebra of $\hat {\frak g}$.  We define its dual $\hat {\frak h}^*$ by 
 \[  \hat {\frak h}^* = {\frak h}^* + {\Bbb C} \delta  + {\Bbb C} L \]
where the pairing of $\hat {\frak h}$ and $\hat {\frak h}^*$ is given by  
 \[ ( a + k_1 K + k_2 d , a' + c_1 \delta + c_2 L ) = ( a , a' ) + k_1 c_2 + k_2 c_1 .         \] 
Let $\Delta_0$ be the set of roots of $\frak g$, the 
 the set of roots of $\hat g$ is
\[ \Delta = \{ \alpha + n \delta \, | \, \alpha \in \Delta_0 , n \in {\Bbb Z} \} \sqcup \{ n \delta \, | \, n\in {\Bbb Z} - \{ 0 \} \} . \] 
 Let $\alpha_1 , \dots , \alpha_n$ be simple roots for $\frak g$;  then
 $\alpha_0 = \delta -\theta ,  \alpha_1 , \dots , \alpha_n$ is a basis of simple roots  for  $\hat {\frak  g}$.

 We apply the above notations to the case $\frak g = sp_{2n}$. 
Let $\rho_0 \in {\frak h}^*$ be such that $( \rho_0 , {\check \alpha}_i )= 1$ for $i=1 , \dots , n$, then 
  $\rho = \rho_0 + (n+1) L \in \hat{\frak h}^*$ satisfies    
$( \rho , {\check \alpha}_i )= 1$ for $i=0, 1 , \dots , n$. 
For each affine real root $ \alpha + n \delta $, its  corresponding 
 root subgroup in $\widehat{\Sp} ( F_v (( t)))$ is  the subgroup $x_{\alpha } ( c t^n ) $ ( $c\in F_v $).
For each element $\lambda \in \hat {\frak h}^*_{\Bbb C} $ of the form 
 \[ \lambda = \lambda_0 + s L + k \delta \]
 where $\lambda_0\in {\frak h}^*$ , 
 we define a quasi-character $\chi_{\lambda}$  of 
 $\sigma ( {\bf A}^* t ) {\rm T}_{\bf A}$ as follows.  For $c_1 , c_2 , c_3 \in {\bf A}^* $, 
\[ 
 \chi_{\lambda } :  c_1 h_{\alpha } ( c_2 ) \sigma  (c_3 t ) \mapsto   (c_1 h_{\alpha } ( c_2 ) \sigma  (c_3 t  ))^{\lambda }
 =  |c_1 |^{s} |c_2|^{ ( \check \alpha , \lambda_0 )  }   | c_3|^ k .\]
  We introduce height functions on 
 $ \widehat {\rm Sp}_{2n} ( {\bf A} \la t \ra ) $ as follows.  For every $g\in \widehat {\Sp}_{2n} ( {\bf A} \la t \ra  )$,
  first we write \[ g = b_g k_g \] by the Iwasawa decomposition,
 with $b_g \in B_{\bf A}, k_g \in K_{\bf A} $, where $B_{\bf A}$ is the restricted product of local Borel subgroups, and
  $ K_{\bf A}$ is the product of local ``maximal compact" subgroups; then we write 
  $   b_g = a_g u_g $ for $a_g \in {\rm T}_{\bf A}$, $u_g $  in the  pro-unipotent radical of $ B_{\bf A}$; 
 then we set 
  \[  \hat h ( \sigma  ( q t  ) g ) = ( \sigma  ( q t  ) a_g )^{ L } \]
 and we put 
\[ \hat h_s ( \sigma  ( q t  ) g ) = ( \sigma  ( q t  ) a_g )^{ s L } =   \hat h ( \sigma  ( q t  ) g )^s  .\]
 It is easy to see that $ \hat h_s( r g k ) = \hat h_s ( g) $ for $r\in \Sp_{2n} ( F \la  t \ra_+  )$ and $ k \in K_{\bf A}$ .

 We consider the Eisenstein series 
\begin{equation}\label{e3} \hat E_s ( T g ) = \sum_{ r \in  \Sp_{2n} ( F\la t \ra_+ ) \backslash   \Sp_{2n} ( F \la t \ra  )} \hat h_s ( r T g ) ,
\end{equation}
  where $T\in {\rm Aut} {\bf A} \la t \ra $ and $g \in \widehat{\Sp}_{2n} ( {\bf A } \la t \ra )$.
  The main result of this section is the following:  

\begin{theorem}\label{convE}
 For $g\in \widehat{\Sp}_{2n} ( {\bf A } \la t \ra ) $ and $ T = ( q_v t )_v \in {\rm Aut} {\bf A} \la t \ra$ satisfying Assumption 4.1,   
 $\hat E_s (T  g ) $ converges absolutely when ${\rm Re} \, s > 3n + 1  $.    
\end{theorem}

The condition  ${\rm Re} \, s > 3n +1 $ is the affine analog of the classical Godement condition for the convergence of 
 Eisenstein series. We use the method in \cite{Garland1} \cite{Garland2} to prove the theorem.  It is enough to prove the 
 theorem for the case $s\in {\Bbb R}$ and $g \in {\rm T}_{\bf A}$, then all the terms in (\ref{e3}) are positive.
 The proof is divided into two steps.

  We first recall that the Bruhat decomposition.  The affine Weyl group $\hat W $ is isomorphic to the semi-direct product
  $ Q^{\vee} \ltimes  W$, where $W$ is the Weyl group of ${\rm Sp}_{2n}$ and $Q^{\vee }$ is the coroot lattice.
We have 
\[   \Sp_{2n} ( F \la  t \ra ) = \sqcup_{w\in W}  B_F  w U_F ,\]
 which implies that 
\[   \Sp_{2n} ( F\la t \ra ) = \sqcup_{w\in Q^{\vee } } \Sp_{2n} ( F \la t \ra_+  ) w U_F , \] 
  where $U$ denote the unipotent radical of $\Sp_{2n} ( F \la t \ra_+  )$.
   Let 
\[  S = \{ w \in W \, | \, w^{-1} \Delta_{0 , + } \subset \Delta_{+ } \} ,\]   
  where $\Delta_{0, +}$ (resp. $\Delta_+$)  is the set of positive roots in $\Delta_0$ (resp. $\Delta $).  It is known that  
$S$ is a set of representatives of  coset space $ W \backslash \hat W$. 
The above decomposition can also be written as 
\[   \Sp_{2n} ( F \la  t \ra ) = \sqcup_{w\in S} \Sp_{2n} ( F \la t \ra_+ ) w U_F . \] 

 Since we assume $s \in {\Bbb R}$, the infinite series (\ref{e3}) has positive terms, the sum makes 
 sense as a function with values in ${\Bbb R}_{ > 0 } \cup \{ \infty \}$.  The group $U $ is a projective limit
  of finite dimensional unipotent groups, and $ U_{F} \backslash U_{\bf A}$ 
 is a projective limit of spaces topologically isomorphic to a finite product of spaces $ F \backslash {\bf A} $. As such
 $ U_{F} \backslash U_{\bf A}$ is a compact topological space, and it 
 has a probability measure invariant under translations by $ U_{\bf A}$. For a given $ g \in \widehat{\rm Sp}_{2n} ( {\bf A} \la t \ra ) $,
 $\hat E_s ( u T  g )$ is an ${\Bbb R}_{\geq 0 } \cup \{ \infty \} $-valued measurable function 
 of $u \in  U_{F} \backslash U_{\bf A}$, because each term is  a measurable function of $u$. 

We first compute for $a\in {\rm T}_{\bf A } \cap \widehat{\rm Sp}_{2n} ( {\bf A} \la t \ra )$,  the constant term
\begin{equation}\label{e4}
  \hat E^{\#}( s,  T a ) =   \int_{  U_F \backslash U_{\bf A} } \hat  E_s ( u T  a ) d u .\end{equation} 
Using the Bruhat decomposition, we have 
\begin{eqnarray}\label{e5}
 && \hat  E_s  ( u T  a )   \nonumber \\
 && = \sum_{w\in S} \, \, \sum_{ r \in     \Sp_{2n} ( F [[t]] ) \backslash \Sp_{2n} ( F[[ t ]] w U_F }      \hat h ( r u \sigma ( qt  ) a )^s 
              \nonumber \\
 && = \sum_{w\in S} \, \, \sum_{ r \in  w^{-1} \Sp_{2n} ( F[[t]]) w \cap U_F \backslash U_F } \hat h ( w r u \sigma ( q t  ) a )^s 
\end{eqnarray}

\noindent Let $\hat E_w ( s ,   u T a ) $ denote the inner sum of (\ref{e5}), and  $\hat E^{\#}_w( s ,  \sigma  ( q t )  a )$
  denote the 
  constant term of $ \hat E_w ( s ,  T a ) $. We have 
 \begin{eqnarray}\label{5aa}
   &&\hat E^{\#}( s,  T  a )    \nonumber \\
    && =   \sum_{w\in S}   E^{\#}_w( s , T  a ) \nonumber \\                      
&&  = \sum_{w\in S}  \, \,  \sum_{ r \in  w^{-1} \Sp_{2n} ( F[[t]]) w \cap U_F \backslash U_F  } \int_{ U_F \backslash U_{\bf A}}     
  \hat h ( w r u T  a )^s    du \nonumber \\
&& =\sum_{w\in S}  \int_{   w^{-1} \Sp_{2n} ( F[[t]]) w \cap U_F \backslash U_{\bf A}  } \hat h ( wu T a )^s    du 
   \end{eqnarray}

 \noindent For each $w\in S$, we put 
 \[ \Delta_{w} =  \{ \alpha \in \Delta_+ \, |\, w \alpha \in \Delta_- \},
 \, \, \, \, \,  \Delta_{w }' =  \{ \alpha \in \Delta_+ \, |\, w \alpha \in \Delta_+ \} . \]
 It is well-known that 
 \begin{equation}\label{sumroot}
 \rho - w^{-1} \rho = \sum_{\alpha \in \Delta_w } \alpha 
\end{equation}

 Let $U_w$ be the group generated by the root subgroups with the roots in $ \Delta_{w}$ and 
 $U_w'$ be the $t$-adic completion  of the group generated by the  root subgroups with the roots in $ \Delta_{w}'$.
 We have the unique factorization 
 $ U = U_w' U_w$ and $ w  U_w' w^{-1} \subset  U $ and 
\[   w^{-1} \Sp_{2n} ( F[[t]]) w \cap U_F =  U_w' ( F ). \]
Therefore 
 \[   w^{-1} \Sp_{2n} ( F[[t]]) w \cap U_F \backslash U _{\bf A}  =       U_w' ( F ) \backslash {U_w}' ({\bf A} )   \cdot  {U_w} ( {\bf A} ).\] 
 We have 
 \begin{equation}\label{e6}
 E^{\#}_w( s , T a ) = \int_{  {U_w} ( {\bf A} )} \hat h  ( w u T a )^s  d u . \end{equation} 
 The right hand side of (\ref{e6}) is 
\begin{eqnarray}\label{e7}
   && \int_{  {U_w} ( {\bf A} )}  \hat h (  w  T a   Ad (\sigma  ( q^{-1} t  ) a^{-1} ) u    )^s d u \\
   &&=   ( \sigma ( qt ) a )^{ w^{-1}( s L ) }  \int_{  {U_w} ( {\bf A} )}  \hat h  (  w    Ad (\sigma  ( q^{-1} t  ) a^{-1} ) u    )^s d u .\nonumber
\end{eqnarray}
Since $   \sigma ( qt ) a   \in \sigma ({\bf A^*} t) {\rm T}_{\bf A} $ normalizes  $  {U_w} ( {\bf A} )$, we make a change
 of variable $ Ad (\sigma  ( q^{-1} t  ) a^{-1} ) u  \to u $, we have 
the right hand side of (\ref{e7}) equals 
 \begin{equation} \label{e8}
   ( \sigma ( q t  ) a  )^{ w^{-1} ( s L ) +  \rho - w^{-1} \rho }    \int_{  {U_w} ( {\bf A} )}  \hat h  ( w u   )^s d u
   \end{equation}

\begin{lemma}\label{cconv} The right hand side of (\ref{5aa}) converges when $ s > 3 n + 1 $ and equals 
 \begin{equation}\label{constant}
E^{\#}( s , T  a ) = \sum_{ w\in S} ( T a  )^{ w^{-1} ( s L ) +  \rho - w^{-1} \rho } 
 \Pi_{\beta  \in  \Delta_+ \cap w  \Delta_-}
 \frac { \hat {\zeta } ( (sL -\rho , \check \beta ) ) } 
 { \hat {\zeta} ( (sL -\rho , \check \beta ) + 1 ) } ,
\end{equation}
where $\hat {\zeta} (s )$ denotes the complete Dedekind zeta function of $F$.
\end{lemma}

\noindent {\it Proof.}  Let $P$ be a finite set of places.  We first consider the integral 
   \begin{equation}\label{e9a}
 I(P)=  \int_{  \Pi_{v\in P} {U_w} ( F_v )}  \hat h  ( w  u   )^s d u  . 
   \end{equation}
  For each $\beta \in \Delta_+ \cap w \Delta_-$, we write $\check \beta = l K + \check \alpha $, 
  where $\check \alpha $ is the coroot of a root of the Lie algebra of $Sp_{2n}$.  Because $ w^{-1} \Delta_{0, +} \subset \Delta_+$, 
  we have $ l\geq 1 $. That implies that 
\begin{eqnarray}
    ( s L - \rho , \check \beta ) &=& ( sL- \rho , l K + \check \alpha ) \nonumber \\
      &=& ( s L , l K + \check \alpha ) - ( \rho , l K + \check \alpha )  \nonumber \\
      &=& l s -(n+1 ) l - (\rho_0 , \check \alpha ) \nonumber \\
      &\geq &  l s -(n+ 1) l - (2 n-1 ) \nonumber \\
      & > &  (3n+1 ) l  - (n+ 1) l - (2 n-1 ) \geq 1 .
 \end{eqnarray}

 Using the Gindikin-Karpelevich formula as in  \cite{Lag}, \cite{Garland1}, we have 
\begin{equation}\label{e9}
  \int_{  {\Pi_{v\in P}  U_w} ( F_v )}  \hat h  ( w  u   )^s d u = \Pi_{\beta }  \Pi_{ v \in P}
 \frac {  {\zeta }_v ( (sL -\rho , \check \beta ) ) } 
 {  {\zeta}_v ( (sL -\rho , \check \beta ) + 1 ) } , \end{equation}
where $\beta $ runs through all  the roots in $ \Delta_+ \cap w  \Delta_-$,  $ \check \beta $ is the coroot corresponding 
 to $\beta $,  and  $ {\zeta }_v ( s ) $ is the $v$-component of the completed Dedekind zeta function of $F$.  
 Let $P$ go to the full set of the places of $F$, we see that 
  the integral 
\[    \int_{ {U_w} ( {\bf A} )}  \hat h  ( w  u   )^s d u \]
 converges to 
\[  \Pi_{\beta  \in  \Delta_+ \cap w  \Delta_-}
 \frac { \hat {\zeta } ( (sL -\rho , \check \beta ) ) } 
 { \hat {\zeta} ( (sL -\rho , \check \beta ) + 1 ) } .\]
    The condition $ |q | > 1 $ implies that 
  \[  ( \sigma ( q t  ) a  )^{ w^{-1} ( s L ) +  \rho - w^{-1} \rho } \]
    decreases to $0$ as fast as a Gaussian function $ e^{ -\epsilon ( v , v ) }$ goes to $0$ as $v \in Q^{\vee }$ gets large,
  so the summation in (\ref{constant}) converges (see \cite{Garland1} for an analogous situation with more details). 
\hfill $\Box$

   After proving the convergence of $E^{\#}( s , T  a )$,
  the 2nd step in the proof of Theorem \ref{convE} is parallel to that of \cite{Garland2}.

 \

 Next we show Theorem \ref{convE} implies the convergence of the Eisenstein series (\ref{Eisen2}):

\begin{theorem}\label{conv2}
   If $m = {\rm dim} \, V > 6n +2 $, then the Eisenstein series (\ref{Eisen2}) converges, consequently 
  the Eisenstein series $ E ( \phi  , Tg  )$ converges for $\phi $ satisfying the conditions in Lemma \ref{compare2} and 
 $T$ satisfying Assumption 4.1.
\end{theorem}

We need to compare the Eisenstein series for two different groups
 $\widetilde {\rm Sp}_{2n} ( {\bf A} \la t \ra )$ and $\widehat {\rm Sp}_{2n} ( {\bf A} \la t \ra )$. 
 We define a third group $ \overline {\rm Sp}_{2n} ( {\bf A} \la t \ra )$ which uses the symbol 
 $ ( a , b ) = | C ( a , b ) |^{- \frac 1 2 } $.
   It has the standard generators $y_{\alpha } ( a ) $ ($\alpha \in \Delta_0$, $a\in {\bf A} \la t \ra ) $).
We have obvious morphisms
 \[ \pi_1 :    \widetilde {\rm Sp}_{2n} ( {\bf A} \la t \ra ) \to   \overline {\rm Sp}_{2n} ( {\bf A} \la t \ra )\]
 given by $ x_{\alpha } (a)  \mapsto y_{\alpha } (a ) ,   c \in {\Bbb C}^* \mapsto  | c |  $,
and 
\[ \pi_2 :    \widehat {\rm Sp}_{2n} ( {\bf A} \la t \ra ) \to   \overline {\rm Sp}_{2n} ( {\bf A} \la t \ra )\]
 given by $ x_{\alpha } (a)  \mapsto y_{\alpha } (a ) ,   c \in {\bf A}^* \mapsto  | c |^{\frac 12 } $.
 We see that 
 for $ r \in  {\rm Sp}_{2n} (  F \la t \ra )$,  $ \pi_1 ( r ) = \pi_2 ( r ) $ and that if  
  $ \pi_1 (g_1 ) = \pi_2 ( g_2 )$, we have 
  $\hat h ( rT g_1 ) = \tilde h ( r T g_1 )^2  $. 
 Therefore $\hat E_s ( T g_1  ) = \tilde E_{ 2s } ( T g_2 )$.  
The condition $ m > 6n +2 $ implies that $\frac m 2  > 3n +1$, by   
 Theorem \ref{convE},   $ \tilde E_{ m  } ( T g_2 ) = \hat E_{\frac m 2} ( T g_1  )$
 is convergent. 

\

\

%
%
 \section{  Eisenstein Series and Snt-modules }

In this section, we show that the $t$-Eisenstein series for snt-modules studied in \cite{gz} appears 
 naturally when we decompose the Eisenstein series $E(\phi , T  g) $ in (\ref{Eisen}). Also we 
 prove Theorem 3.3 in \cite{gz} which was stated there without proof.   

  For $T$ as in Assumption 4.1,  $ \phi $ as in Lemma \ref{compare2},
 we have the Eisenstein series $E ( \phi , T g )$ which is convergent by Theorem \ref{conv2}.
 For simplicity we  put $f = T  g \phi $.
 We have 
\begin{equation}\label{6.1}
 E (  T g \phi  )  =  E( f ) = \sum_{ r\in  {\rm Sp}_{2n} ( F\la t \ra_+ ) \backslash    {\rm Sp}_{2n} ( F\la t \ra ) }
     ( r f) ( 0 )        .
 \end{equation}

\noindent Let \[        Gr ( F(( t))^{2n} ) \]
denote the set of Lagrangian subspaces $U$ satisfying the following conditions (1). $U$ is an $F[[t]]$-submodule; (2). $U$ is commensurable 
  with $ F[[t]]^{2n}$. 

\begin{lemma}\label{lemma6.0} ${\rm Sp}_{2n } ( F  ((  t )) )$ acts on $Gr ( F(( t))^{2n}$ transitively, the isotropy subgroup of 
  $F[[t]]^{2n}$ is $ {\rm Sp}_{2n} ( F [[ t ]]  )$.
\end{lemma}

This lemma is well-known. Since we could not find the reference, we give a sketch of the proof.  

\noindent{\it Proof.}  For $ U \in  Gr ( F(( t))^{2n} )$. We first note that for every $a, b \in U$, 
 the $F((t))$-valued sympletic paring $\la a , b \ra_{ F((t))}$ is in $F[[t]]$. Otherwise 
 there is a maximal, positive integer $k$ such that the $t^{-k}$-coefficient of   $\la a , b \ra_{ F((t))}$ is  not $0$; then 
  $ \la a , t^{k-1} b \ra = {\rm Res} \la  a , t^{k-1} b \ra_{ F((t))} \ne 0 $, which contradicts  $a$ and $t^{k-1}b$ being elements in 
  the Lagrangian subspace $U$.  Since $U$ is commensurable with  $ F[[t]]^{2n}$, it is isomrophic to $ F[[t]]^{2n}$ as an $F[[t]]$-module.
 We can find an $F[[t]]$-basis $a_1 , \dots , a_{2n}$ such that $ \la  a_i ,  a_{j +n } \ra_{ F((t))} = \delta_{i, j}$ for 
  $1\leq  i , j \leq n $.  Let $e_1 , \dots , e_{2n}$ be the standard sympletic basis of the $F((t))$-symplectic space $F((t))^{2n}$, 
 then the $F((t))$-linear map $g:  e_i \mapsto a_i $ ( $ 1 \leq i \leq 2n$) is in ${\rm Sp}_{2n } ( F  ((  t )) )$ and 
   $ F[[t]]^{2n}  g  = U$.  It is clear that the isotropy subgroup of 
  $F[[t]]^{2n}$ is $ {\rm Sp}_{2n} ( F [[ t ]]  )$.    \hfill $\Box$

  We next prove that the above lemma is also true if  ${\rm Sp}_{2n } ( F(( t)))$ and ${\rm Sp}_{2n} ( F[[t]] )$
 are replaced by smaller groups ${\rm Sp}_{2n } ( F \la  t \ra )$ and ${\rm Sp}_{2n} ( F\la t \ra_+  )$.

\begin{lemma}\label{lemma6.1} ${\rm Sp}_{2n } ( F \la  t \ra )$ acts on $Gr ( F(( t))^{2n}$ transitively, the isotropy subgroup of 
  $F[[t]]^{2n}$ is $ {\rm Sp}_{2n} ( F\la t \ra_+  )$.
\end{lemma}

 \noindent{\it Proof.} Using the BN-pair  argument as in \cite{S}, we have 
 \[ {\rm Sp}_{2n } ( F ((  t  ) )) = {\rm Sp}_{2n } ( F [[  t ]] )
 {\rm Sp}_{2n } ( F (  t , t^{-1}) ). \] 
It is clear that $  {\rm Sp}_{2n } ( F (  t , t^{-1}) ) \subset  {\rm Sp}_{2n } ( F \la t \ra )$.
 The transitivity follows. 
  Since the isotropy subgroup of the Lagrangian subspace $F[[t]]^{2n}$ in $ {\rm Sp}_{2n } ( F ((  t  ) ))$ is
  ${\rm Sp}_{2n } ( F [[  t ]] )$, its isotropy subgroup 
 in ${\rm Sp}_{2n } ( F \la t \ra )$ is 
 \[ {\rm Sp}_{2n } ( F [[  t ]] ) \cap {\rm Sp}_{2n } ( F \la t \ra ) = {\rm Sp}_{2n } ( F \la t \ra_+ ) .\]
\hfill $\Box$

\

 By Lemma \ref{lemma6.1}, the summation in (\ref{6.1}) can be written as a summation over $   Gr ( F(( t))^{2n} ) $.  For 
$ r \in  {\rm Sp}_{2n} ( F\la t \ra_+ ) \backslash    {\rm Sp}_{2n} ( F\la t \ra )$, let $ U = F[[ t]]^{2n} r $ be its 
 corresponding element in $Gr ( F(( t))^{2n} ) $.  Let $\pi_- :  F((t))^{2n} \to  t^{-1} F [ t^{-1} ]^{2n}$ be the projection
  map with respect to the decomposition:
\begin{equation}\label{6.1.1}  F((t))^{2n}  =  F[[t]]^{2n} +  t^{-1} F [ t^{-1} ]^{2n}. \end{equation}
 By Proposition \ref{prop}, we have  
\begin{equation}\label{6.2}
  ( r f  )   (0 ) = \int_{ {\rm Im } (\gamma_r)_{\bf A}}  \psi ( \frac 12 \la x^* \gamma_r , x^* \delta_r \ra ) 
      f ( x^* \gamma_r ) d ( x^* \gamma_r ) .
\end{equation}
The symplectic pairing 
\[ \la , \ra :  \pi_- ( U ) \times F[[t]]^{2n} \to F \] 
factors through a non-degenerate pairing 
\[ \la , \ra :  \pi_- ( U ) \times F[[t]]^{2n}  /  F[[t]]^{2n} \cap U      \to F .\]
For each $ v \in \pi_- ( U )$, let $\tilde v \in U $ be a lifting of $v$, write $\tilde v = \tilde v_+  + \tilde v_- $ according to 
 the decomposition (\ref{6.1.1}),
 then the element $  \tilde v_+ +   F[[t]]^{2n} \cap U  \in F[[t]]^{2n}  /  F[[t]]^{2n} \cap U  $ is independent of the lifting.
 We denote  by $\rho $ the map: 
\[ \rho:   \pi_- ( U ) \to  F[[t]]^{2n}  /  F[[t]]^{2n} \cap U ,   \, \, \,  v \mapsto \tilde v_+ +   F[[t]]^{2n} \cap U .\]
And we use the same symbol $\rho$ to denote the map 
\[ \rho = \rho \otimes Id :    \pi_- ( U )\otimes V  \to  (F[[t]]^{2n}  /  F[[t]]^{2n} \cap U)\otimes V . \]
Then (\ref{6.2}) can be written as  
\begin{equation}\label{eu}
  ( r  f)  (  0 ) = \int_{ (\pi_- ( U ) \otimes V )_{\bf A}} \psi (  \frac 12 \la x , \rho x \ra ) f ( x ) d x \deff E( f , U )  ,
\end{equation}
where $dx$ is the Haar measure on  $(\pi_- ( U ) \otimes V)_{\bf A}$ such that the covolume of $\pi_- ( U ) \otimes V$ is $1$.
We have 
\[  E ( T g \phi ) = E( f ) = \sum_{ U \in Gr ( F((t))^{2n} ) } E ( f , U ) . \]

We consider the space $ t^{-1} F [ t^{-1} ]^{2n} $ as an $F[[t]]$-module by the identification
 $t^{-1} F [ t^{-1} ]^{2n} =  F (( t))^{2n} / F [[ t]]^{2n}$. 
The group ${\rm Sp}_{2n} ( F[[t]]) $ acts
 on  $ t^{-1} F [ t^{-1} ]^{2n} $ by means of the above identification.

We denote by
 \[ Gr (  t^{-1} F [ t^{-1} ]^{2n} ) \]
 the set of all $F[[t]]$-submodules in $t^{-1} F [ t^{-1} ]^{2n}$ which are  finite dimensional as an $F$-space.
  The projection map $\pi_- : F(( t))^{2n} \to t^{-1} F [ t^{-1} ]^{2n}$ with respect to the decomposition (\ref{6.1.1})
 gives  a map 
\[  P :  Gr ( F(( t))^{2n} ) \to Gr (  t^{-1} F [ t^{-1} ]^{2n} ) , \, \,    U \mapsto \pi_- ( U ) .\]
The map $P$ is ${\rm Sp}_{2n} ( F [[ t]])$-equivariant  but is not surjective. For example, 
 $  F e_1 t^{-1} + F e_{n+1 } t^{-1} \in  Gr (  t^{-1} F [ t^{-1} ]^{2n} ) $, but is not in the image of $P$. Otherwise, if
  $F e_1 t^{-1} + F e_{n+1 } t^{-1} = \pi_- ( U) $,  then $U$ contains elements
 $ e_1 t^{-1} + a  , e_{n+1 } t^{-1} + b $ for some $a , b \in F[[t]]^{2n}$. Since $U$ is an $F[[t]]$-module, 
  $ e_{n+1} + t b \in U$,  and then 
\[ \la e_1 t^{-1} + a , e_{n+1} + t b \ra = \la e_1 t^{-1}, e_{n+1 } \ra = 1 , \]
 which contradicts  $U$ being Lagrangian.

For $W\in Gr ( t^{-1} F [t^{-1}]^{2n} )$, we set
 \[E_W ( f) =     \sum_{ U \in Gr ( F((t))^{2n} ) : P( U ) = W }  E ( f , U) .\]
  In the case that $ P^{-1} ( W)$ is empty, then we set $E_W ( f)  = 0$.  
We have 
 \begin{equation}\label{6de}   E ( f ) = \sum_{ W \in Gr ( t^{-1} F [t^{-1}]^{2n} )}  E_W ( f) .
 \end{equation}

  The following lemma describes the image of $P$.

\begin{lemma}\label{lemma6.2}  An element $ W \in Gr (  t^{-1} F [ t^{-1} ]^{2n} ) $ is in the image of the map $P$ iff
 there is $g\in {\rm Sp}_{2n} ( F [[ t]] )$ such that 
 \begin{equation}\label{6.2c}
  W g = {\rm Span}_{ F[[t]] } (  t^{-k_1 } e_1 ,  t^{-k_2} e_2 ,  \dots  , t^{-k_l } e_l ) \end{equation}
where $l \leq n$ and  $k_1 \geq k_2 \geq \dots \geq k_k \geq 1 $.
\end{lemma}

\noindent {\it Proof.}   Let $U$ denote the $F[[t]]$-submodule of $F((t))^{2n}$ generated by 
   \[   t^{-k_1 } e_1 ,  t^{-k_2} e_2 , \dots  ,t^{-k_l } e_l , e_{l+1} , \dots ,e_n , t^{ k_1  } e_{n+1} , \dots,
      t^{k_l} e_{n+l} , e_{n+l+1} , \dots , e_{2n} .\]
 It is easy to see that $U \in Gr ( F(( t))^{2n} ) $. Then 
   \[ P ( U ) = {\rm Span}_{F[[t]]} (  t^{-k_1 } e_1 ,  t^{-k_2} e_2 , \dots  ,t^{-k_l } e_l ) . \]
 This and the ${\rm Sp}_{2n} ( F[[t]])$-equivariance of $P$ prove the condition in the Lemma is sufficient.
  Conversely if $ W = P(U)$ for $U \in Gr ( F((t))^{2n})$, then $ U = F[[t]]^{2n} g $ for some $g\in  {\rm Sp}_{2n} ( F((t)))$.
  By the Bruhat decomposition,  we write $g$ as 
  \[ g = b_1 \, {\rm diag} ( t^{-k_1 } , \dots , t^{-k_n } , t^{k_1 } \dots , t^{k_n} ) \, b_2 \]
 for some $ k_1 \geq k_2 \geq \dots \geq k_n \geq 0 $, $b_1 , b_2 \in {\rm Sp}_{2n} ( F[[t]])$.  Let $k_l$ be the 1st in $k_i$'s that is not $0$, then 
  \[ P( U ) =   {\rm Span}_{ F[[t]] } (  t^{-k_1 } e_1 , t^{-k_2} e_2 , \dots  , t^{-k_l } e_l ) \cdot  b_2 \]
This proves the condition is also necessary. \hfill $\Box $

\

For $W\in Gr (  t^{-1} F[t^{-1} ]^{2n} ) $,   $W + F[[t]]^{2n}$ is an $F[[t]]$-submodule of $F((t))^{2n}$. 
  Since the symplectic form $\la , \ra$ on $F((t))^{2n}$ satisfies the property $ \la t a , b \ra = \la a , t b \ra $,
the radical $R$ of  the restriction of
 $\la , \ra $ on $W + F[[t]]^{2n}$ is an $F[[t]]$-submodule and $ R \subset F[[t]]^{2n} $. 
The quotient $F[[t]]$-module 
 \[  \tilde W \deff  ( W + F[[t]]^{2n} ) / R \]
has the induced  symplectic form $\la , \ra $    
 and it satisfies the condition $ \la t a , b \ra =  \la  a , t  b \ra$,
 i.e.,   $\tilde W$ has   the structure of an snt-module (See Section 1 \cite{gz}  for definition of an snt-module).
  In the case that $ W$ is the right hand side of  (\ref{6.2c}), the snt-module  $\tilde W$ is isomorphic to 
   \[  H_{k_1 } \oplus \dots \oplus  H_{k_l } ,\]
 where $H_k $ is as in (1.9) \cite{gz}.
Recall from \cite{gz}, Section 2, that for an snt-module $\tilde W$, we used the symbol $ Gr ( \tilde W , t)$ to denote the   
 set of Lagrangian subspaces of $\tilde W$ which are also $F[[t]]$-submodules.  We now prove that there is a bijection
 from    $ Gr ( \tilde W , t)$ to the set  
 
\begin{equation}\label{6.4} 
  Gr ( F((t))^{2n} )_{\leq W} \deff   \{  U \in Gr ( F((t))^{2n} )\, \, | \, \, \pi_- ( U ) \subset W \}.
\end{equation} 

 \
   
\noindent  For $  U $ as in (\ref{6.4}), since  $\pi_- ( U) \subset W$, we have  $ U \subset  W + F[[t]]^{2n}$,  
 so $ U /  U \cap R $ is an $F[[t]]$-submodule of $ \tilde W $. And the fact that $U$ is a Lagrangian subspace of $F(( t))^{2n}$
 implies that $U /  U \cap R $ is a Lagrangian subspace of $\tilde W$. Therefore $  U /  U \cap R \in Gr ( \tilde W  , t )$.
 Conversely if $ M \subset \tilde W$ is an element in $ Gr ( \tilde W , t)$,  then 
  the inverse image of $M$ under the canonical map $  W + F[[t]] \to \tilde W$ is an element in $Gr ( F((t))^{2n} )$. 
 This proves 

\begin{lemma}\label{lemma6.3} The map $  Gr ( F((t))^{2n} )_{\leq W} \to Gr ( \tilde W , t )$ given by 
  $U \mapsto  U / U \cap R $ is a bijection.
\end{lemma}

\noindent  In the snt-module $\tilde W$, $ F[[t]]^{2n} / R \in Gr ( \tilde W , t )$, and 
    $W / W \cap R = W$ is a Lagrangian subspace of $\tilde W$, but in general $ W $ is not an $F[[t]]$-submodule of $\tilde W$,
 i.e., $ W \notin Gr ( \tilde W , t )$.  We have decomposition 
  \[   \tilde W  =  W \oplus  F[[t]]^{2n} / R \]
 into a sum of Lagrangian subspaces.  As in Section 8 \cite{gz} we have 
  \[ {\cal S} ( ( W \otimes V)_{\bf A}) \]
 as a model of the Weil representation for the symplectic group ${\rm Sp} ( (\tilde W \otimes V)_{\bf A} )$, 
 and the groups $ {\rm Sp} ( \tilde W , t )_{\bf A}$ and $ G^q ( F[[ t]] ) _{\bf A})$ form a commuting pair
 in ${\rm Sp} ( (\tilde W \otimes V)_{\bf A} )$. 
 For each $ f\in  {\cal S} ( ( W \otimes V)_{\bf A}) $, we defined in \cite{gz} the Eisenstein series
 
 \[  {\rm Et}  ( f ) = \sum_{ H \in  Gr ( \tilde W , t ) } E ( f , H ).  \]

Compare the formula $ E ( f , H) $ in Section 8  \cite{gz} and the formula $ E ( f , U)$ in (\ref{eu}) above,
 we see that if $H$ corresponds to $U$ in the correspondence in Lemma \ref{lemma6.3}, then 
 $ E ( f , H) = E( f , U)$.

 Therefore we have

\begin{lemma}\label{lemma6.4} We use the same symbol $f$ to denote the restriction of $f$ on $ (W \otimes V)_{\bf A}$, then 
  $t$-Eisenstein series $ {\rm Et} ( f )$ for the snt-module $\tilde W$ is 
 \[    {\rm Et} ( f ) = \sum_{ U \in Gr ( F(( t))^{2n})_{\leq W } } E ( f , U ). \]
 And ${\rm Et}_W ( f ) $ in Section 8 \cite{gz} is 
\[{\rm Et}_W ( f ) = E_W ( f )  .\]
\end{lemma}

  When $ f \in {\cal S}  ( (W \otimes V )_{\bf A})$ is the restriction of $ T g \phi \in  {\cal S}  ( ( t^{-1} F[t^{-1}]^{2n}\otimes V )_{\bf A})$
for $ T , g , \phi $ as in Theorem \ref{conv2} , since ${\rm Et} ( f )$ is part of series 
 for $ E ( f )$, we know from Theorem \ref{conv2} that $E( f )$ converges absolutely, therefore ${\rm Et} ( f )$ converges
 absolutely. We shall prove a stronger convergence result:

\begin{theorem} \label{thm6.1}  Suppose ${\rm dim} V > 6 n +2 $ and 
 $W \in {\rm Im } P $,  then 
 $ {\rm Et} ( \phi )$ converges absolutely and the convergence is uniform for $\phi$ in a compact subset of ${\cal S} ( (W\otimes V)_{\bf A} )$. 
\end{theorem} 

This Theorem  implies Theorem 3.3  \cite{gz}.  
To prove Theorem \ref{thm6.1}, we introduce some useful terminology. 
 A subset ${\cal C} \subset {\cal S} ( ( t^{-1} F[t^{-1}]^{2n}\otimes V )_{\bf A})$ is called {\it quasi-compact} if 
 it satisfies the following two conditions: (1) there is a 
 finite set $S$ of places of $F$ including all infinite places such that for $v\notin S$, every $f\in {\cal C} $ is of the 
 form $ f = f' f_0 $ where $f'$ is in ${\cal S} ( \Pi_{v\in S} ( t^{-1} F[t^{-1}]^{2n}\otimes V )_{F_v} )$ and
 $f_0 = \Pi_{v\notin S} \phi_{0 , v } $, where $\phi_{0, v}$ is as Lemma \ref{phi0padic}.
   (2) For every finite dimensional space $W \subset   t^{-1} F[t^{-1}]^{2n}\otimes V $,
 the restriction of 
 ${\cal C}$ to $ W_{\bf A}$ is a subset of some compact subset of ${\cal S} ( W_{\bf A} )$. 
 We give an example of a quasi-compact subset.  Let $t^{-1} F[t^{-1}]^{2n}\otimes V =  W' \times W'' $ be a decomposition 
  of vector spaces such that $W'$ is finite dimensional.  Let $ {\cal C}_1 \subset {\cal S} ( W'_{\bf A} )$ be a compact subset,
 and $f_2 \in {\cal S} ( W''_{\bf A})$ be a fixed function, then for each $f_1 \in {\cal C}_1 $, $ f_1 ( x_1 ) f_2 ( x_2 )$ (where
 $x_1 \in W'_{\bf A}, x_2 \in W''_{\bf A} $ ) is a function on $ ( W' \times W'')_{\bf  A} $, so 
${\cal C}_ 1 f_2 $ can be regarded as a subset 
 of $ {\cal  S} (	( W' \times W'' )_{\bf  A} ) ={ \cal S} ( ( t^{-1} F[t^{-1}]^{2n}\otimes V )_{\bf A} )$, this set is    
clearly quasi-compact. 
We also have the concept of quasi-compact subset of $ {\cal S} ( t^{-1} F_v [ t^{-1}]^{2n} \otimes V_v) $: a subset is called quasi-compact 
 if its restriction on each finite dimensional subspace $W$ is a subset of a compact subset in ${\cal S} ( W)$.

 Now we fix a place $v$.  Recall in Section 2, we defined for each simple root $\alpha_i $ ($i=0, 1 , \dots , n$) , a subgroup 
  $K_{\alpha_i } \subset \widetilde {\rm Sp}_{2n} ( F_v ((t)) )$. 
  
\begin{lemma}\label{compact1} If ${\cal C} \in {\cal S} ( t^{-1} F_v [ t^{-1} ] ^{2n} \otimes V_v )$ is a quasi-compact set,
  the for each simple root $\alpha_i$,  $ K_{\alpha_i } {\cal C} $ is also quasi-compact.
\end{lemma}

\noindent {\it Proof.} If $i=1 , \dots , n$, since $ K_{\alpha_i } \subset {\rm Sp}_{2n} ( F_v )$ is a compact subset, and 
  for $ g \in {\rm Sp}_{2n} ( F_v )$, $ f \in {\cal C}$, 
  \[ ( g f ) ( x ) = f ( x g ) .\]
It is clear that $ K_{\alpha_i } {\cal C}$ is quasi-compact.  It remains to prove the Lemma for the case $K_{\alpha_0 }$. 
 In the notation of Section 2.2, $K_{\alpha_0} $ is a maximal compact subgroup $G_{\alpha_0}$, and 
  $G_{\alpha_0}$ is generated by root vectors $ x_{\theta} ( c t^{-1} )$ and $  x_{-\theta} ( c t )$ ( $ c\in F_v $).
  We first note that if $ B$ is  a compact subset in $ F_v$,  then the following subsets
 \[     \{ x_{\theta} ( c t^{-1} ) \, | \, c \in B \}{\cal C} , \, \, \, \, \, \, 
   \{ x_{- \theta} (  c t  ) \, | \, c \in B \}{\cal C} \]
 are quasi-compact.
This can be easily seen from the action formula of $x_{\theta} ( c t^{-1} )$ and $x_{- \theta} (  c t  )$ in (\ref{2a1}).
 We notice that a partial Fourier transform maps a quasi-compact set to a quasi-compact set, and a map that permutes the variables
 maps a quasi-compact set to a quasi-compact set.
 Having noticed these facts, we divide our proof into two cases. If $v$ is a finite place, let ${\cal O}_v$ denote the ring of 
 integers of $F_v$.
Consider each element
\begin{equation}  g = \left( \begin{matrix}  u & v \\ w &  s \end{matrix} \right)  \in SL_2 ( {\cal O}_v ). \end{equation} 
 We claim $g$ is in one of the following two  sets: 
\[   \{ \left( \begin{matrix}  1 & 0 \\ c  &  1 \end{matrix} \right) 
    \left( \begin{matrix}  h & 0 \\ 0 &  h^{-1}  \end{matrix} \right) 
    \left( \begin{matrix}  1 & b \\ 0 &  1 \end{matrix} \right) \, |  \, h \in {\cal O}_v ^* ,  b , c\in  {\cal O}_v  \} \]
\[ \{\left( \begin{matrix}  1 &  1 \\ 0  &  1 \end{matrix} \right) \left( \begin{matrix}  1 & 0 \\ c  &  1 \end{matrix} \right) 
    \left( \begin{matrix}  h & 0 \\ 0 &  h^{-1}  \end{matrix} \right) 
    \left( \begin{matrix}  1 & b \\ 0 &  1 \end{matrix} \right) \, |  \, h \in {\cal O}_v ^* ,  b , c\in  {\cal O}_v  \} \]
In fact, if $ u\in {\cal O}_v ^* $, then $ g $ is in the first set, if $u \notin {\cal O}_v^* $, then $ v ( u ) \geq 1$, 
 since $ u s - v w =1 $, we must have $ w \in {\cal O}_v^*$, then 
   \[  \left( \begin{matrix}  1 &  -1 \\ 0  &  1 \end{matrix} \right) g  \]
is in the first set, so $ g $ is in the second set.  Therefore every $ g \in K_{\alpha_0 }$ is in one of the two sets
\[  B_1 =  S^1 \{  x_{\theta } ( c t^{-1} ) h_{\theta} ( h )  x_{-\theta } ( b t ) \, | \,  h \in {\cal O}_v^* , b , c \in {\cal O}_v \} \]
\[ B_2 =  S^1 \{ x_{-\theta } ( t  )  x_{\theta } ( c t^{-1} ) h_{\theta} ( h )  x_{-\theta } ( b t ) \, | \,  h \in {\cal O}_v^* , b , c \in {\cal O}_v \}\]
It is clear that $B_1 {\cal C}$ and $ B_2 {\cal C}$ are quasi-compact.  So 
 $ K_{ \alpha_0 } {\cal C} \subset  B_1 {\cal C} \cup  B_2 {\cal C}$ is also quasi-compact.
  For the case $F_v={\Bbb R}$ and for every $g \in K_{\alpha_0}$, let   
 \[ \left( \begin{matrix}  {\rm cos} \, \theta  &  {\rm sin} \, \theta  \\ - {\rm sin} \, \theta    &  {\rm cos} \, \theta\end{matrix} \right) \]
 be its  image in $SO_2$.  If $ |   {\rm cos} \, \theta  | \geq \frac {\sqrt 2}  2  $, using the identity
 \[    \left( \begin{matrix}  {\rm cos} \, \theta  &  {\rm sin} \, \theta  \\ - {\rm sin} \, \theta    &  {\rm cos} \, \theta\end{matrix} \right) 
   =    \left( \begin{matrix}  1  &  0  \\ - \frac { {\rm sin} \, \theta } { {\rm cos}\, \theta }
        &   1   \end{matrix} \right) 
     \left( \begin{matrix}  {\rm cos} \, \theta  &  0  \\  0 & { {\rm cos}^{-1} \, \theta }
          \end{matrix} \right) 
 \left( \begin{matrix}  1  &   \frac { {\rm sin} \, \theta }{{\rm cos} \, \theta } 
  \\ 0   &   1   \end{matrix} \right) , \]
we see that $ g$ is in
  \[  B_3 \deff  S^1 \{  x_{\theta } ( c t^{-1} ) h_{\theta} ( h )  x_{-\theta } ( b t ) \, | \,  \frac {\sqrt 2}  2 \leq    |h | \leq 1   , \; | b |\leq 1  ,\; |c | \leq  1 \} .\]
  If  $ |   {\rm cos} \, \theta  | \leq \frac {\sqrt 2}2  $, then $ w_{\alpha_0} ( 1) g $ is in $ B_3$ above. 
  Therefore $ K_{\alpha_0} {\cal C} \subset   B_3 {\cal C} \cup  w_{\alpha_0} (1) B_3 {\cal C}$ which is quasi-compact. 
 A similar proof works for $F_v = {\Bbb C}$.  
\hfill $\Box $

\

\begin{lemma}\label{compact2} Let $W \in {\rm Im } P $, then there exists a positive integer $l$ (depending on $W$) such that 
         each $  r \in Sp_{2n} ( F \la t \ra )$ satisfying the condition that $   P (   F[[t]]^{2n} r ) \subset W $
      can be written as 
  \[ r = b k_{1 } k_{2 } \dots k_{l } \]
 where $ b\in B_{\bf A}$ and $ k_{i} $ is in $  \Pi_v  K_{v, \alpha }$ for some simple root $\alpha$.
\end{lemma}

\noindent {\it Proof.}  
   We may assume that $ W = {\rm Span}_{ F[[ t]]}  \{  t^{-j_1 } e_1 , t^{-j_2 } e_2 , \dots , t^{-j_m } e_m \} $
 where $ m \leq n $.  By the Bruhat decomposition, every $ r $ can be written as 
 $ r = b_1 w b_2 $ for $b_1 , b_2 \in B_F $ and $w \in W$. The condition on $r$ implies that the length $l( w ) $ of $w$ is bounded, say 
  $ l( w ) \leq l $ for all $r$ satisfying the condition in the lemma. Then we use 
   the BN-pair and the argument as in p.99 \cite{S} to show that $r$ can be written as
  \[ r = b k_{1 } k_{2 } \dots k_{l } \] with $ k_i \in \Pi_v K_{v , \alpha_i }$ for some simple root $\alpha $.
\hfill $\Box $

\begin{lemma}\label{lemma6.9} Let ${\cal C}$ be a quasi-compact subset of ${\cal S} ( (t^{-1} F[t^{-1}]^{2n} \otimes V)_{\bf A} ) $,
 let $W \in {\rm Im } P $.  There is a constant $C$ depending only on ${\cal C}$ and $W$ such that 
 for every $r\in {\rm Sp}_{2n} ( F\la t \ra )$ satisfying the condition that $ \pi_- ( F[[t]]^{2n} r ) \subset W$, 
 $f\in {\cal C}$, we have 
   \[ ( r T f ) ( 0 ) \leq C \tilde h ( r T ) ^m .\]
\end{lemma}

\noindent {\it Proof.}  By Lemma \ref{compact2}, we can write $ r = b k_1 \dots k_l $, so 
     \[  ( r T f ) ( 0 ) = \tilde h ( r T )^{m } ( k_1 \dots k_l f ) ( 0 ) .\]
By Lemma \ref{compact1}, $ k_1 \dots k_l f $ varies over some quasi-compact subset ${\cal C}'$; then there is a constant $C$ satisfying   
   \[  |( k_1 \dots k_l f ) ( 0 )| \leq  C \]
for all $f\in {\cal C}$ .  \hfill $\Box $.

\

\noindent {\it Proof of Theorem \ref{thm6.1}.}  Consider ${\cal C}_0 $
 a compact subset of $ {\cal S}  (( W\otimes V )_{\bf A}) $, we choose a complement 
   $W'$ of $W$ in $ t^{-1} F [ t^{-1} ]^{2n} $, i.e., $ W\oplus W' =  t^{-1} F [ t^{-1} ]^{2n}$,  and  we 
  choose a function $ f \in   {\cal S}  (( W'\otimes V )_{\bf A}) $, then $ T^{-1} {\cal C}_0 f $ is a quasi-compact subset of 
  $ {\cal S} ( ( t^{-1} F [ t^{-1} ]^{2n}\otimes V)_{\bf A} ) $.  Consider for $\phi \in {\cal C}_0$,  
  \[ {\rm Et}_W ( \phi ) = \sum_{ r \in {\rm Sp}_{2n} ( F[[ t ]] )  \backslash  {\rm Sp}_{2n} ( F(( t)) ):   P ( F[[t]]^{2n} r ) \subset W } 
    ( r T T^{-1} \phi  ) ( 0 ) \] 
  By Lemma \ref{lemma6.9}, there is a constant $C$ such that 
  \[ | ( r T T^{-1} \phi  ) ( 0 ) | \leq C \tilde h ( r T )^{m }.  \] 
    Therefore 
 \[ | {\rm Et}_W ( \phi ) | \leq C \sum_{ r \in {\rm Sp}_{2n} ( F[[ t ]] )  \backslash  {\rm Sp}_{2n} ( F(( t)) ):   P ( F[[t]]^{2n} r ) \subset W }
    \tilde h ( r T )^{m} ,\]
 since $ m > 6 n +2 $, the right hand side is convergent by Theorem \ref{conv2}. 
\hfill $\Box $

\

\

%
%
%
%

\section{Siegel-Weil Formula for Loop Groups}

In this section, we prove our main result: the Siegel-Weil formula for loop groups.
Let  $ T , \phi $ be as in Theorem \ref{conv2}, and $h\in G( {\bf A}  \la t \ra_+ ) $.  We consider the theta functional 

 \begin{equation}\label{7.1} 
 \theta ( h T g \phi ) =  \sum_{r \in t^{-1} F[ t^{-1} ]^{2n} } ( h T g \phi ) ( r ).
 \end{equation}

\noindent We fix $ T , g , \phi $ and denote $ f = T g \phi $ and 
  regard (\ref{7.1}) as a function of  $ h \in  G( F  \la t \ra_+ )  \backslash G( {\bf A} \la t \ra_+ )$.

\begin{lemma}\label{lemma7.1} The convergence of series (\ref{7.1}) is uniform for $h\in G ( {\bf A} \la t \ra_+ )$.
\end{lemma}

  Note that the assumption that $G$ is the orthogonal group of an anisotropic form implies that
 $G( F \la t \ra_+ )  \backslash G( {\bf A} \la t \ra_+ )$ has a compact fundamental domain.
 The proof of this lemma is similar to that of Theorem \ref{converge}, where Lemma 3.5 is used. 
 We also remark that $G( F \la t \ra_+ )  \backslash G( {\bf A} \la t \ra_+ ) = G( F [[ t ]] )  \backslash G( {\bf A} [[ t ]] )$.

 Since the convergence of (\ref{7.1}) is uniform on $h \in G( {\bf A} \la t \ra_+ )$, this function is continuous.
 And since we assume the quadratic space $V$ is anisotropic over $F$,  $ G( F \la t \ra_+ )  \backslash G( {\bf A} \la t \ra_+ )$
 is compact,  the following  integration

\begin{equation}\label{7.2} 
I ( Tg \phi ) = I ( f ) \deff \int_{ G( F \la t \ra_+ )  \backslash G( {\bf A} \la t \ra_+ ) } \theta ( h f ) d h ,
\end{equation}
converges, 
where $dh $ is the unique $ G( {\bf A} \la t \ra_+ )$-invariant probability measure on  $G( F \la t \ra_+ )  \backslash G( {\bf A} \la t \ra_+ )$.

  We write  $ h \in G ( {\bf A} \la  t \ra_+ )$ in the block form   
 \begin{equation*}
 \left[ 
 \begin{array}{cc}
 \alpha _{h} & \beta _{h} \\ 
 \gamma _{h} & \delta _{h}%
 \end{array}%
 \right]
  \end{equation*}%
with respect to the decomposition 
   \[ (F (( t)) ^{2n} \otimes V)_{\bf A} = (t  F[ t^{-1} ]^{2n} \otimes V)_{\bf A} \oplus ( F[[ t]]^{2n} \otimes V)_{\bf A} .\]
 Since $ G ( {\bf A} [[ t ]] )$ preserves the second summand,   we have 
 $\gamma_h= 0 $.
 Then $h \cdot f $  is given by 
\begin{equation}\label{7.3}
  ( h  \cdot f ) ( r ) =  \psi ( \frac 12 \la r \alpha_h , r \beta_h \ra )  f ( r \alpha_h  ).
\end{equation}
Therefore, (\ref{7.2}) can be written as 
\begin{equation}\label{7.4} 
  I ( f ) = \int_{ G( F \la t \ra_+ )  \backslash G( {\bf A} \la t \ra_+ ) }
 \sum_{r \in t^{-1} F[ t^{-1} ]^{2n} } \psi ( \frac 12 \la r \alpha_h , r \beta_h \ra )  f ( r \alpha_h  )dh.
 \end{equation}

We can write the integral (\ref{7.4}) as a sum of orbital integrals. We first 
define  right actions of $G( {\bf A} \la t \ra_+ )$ on $(t^{-1} F[t^{-1} ]^{2n} \otimes V)_{\bf A} $
and $ S^1 \times  (t^{-1} F[t^{-1} ]^{2n} \otimes V)_{\bf A}   $ as follows:
for $ ( s , x ) \in   S^1 \times  (t^{-1} F[t^{-1} ]^{2n} \otimes V)_{\bf A} $,
  $h \in  G( {\bf A} \la t\ra_+)$
\[   x h = x \alpha_h  , \, \, \, \, \,  ( s, x ) h =  ( s \psi (  \frac 12 \la x \alpha_h , x \beta_h \ra ) ,  x \alpha_h  ) .\]
One checks directly that the above are actions. 
It is clear that the projection map  
  \[ S^1 \times  (t^{-1} F[t^{-1} ]^{2n} \otimes V)_{\bf A} \to (t^{-1} F[t^{-1} ]^{2n} \otimes V)_{\bf A} , \, \, \, \, \, ( s, x ) \mapsto x \]
is $G( {\bf A} [[t]])$-equivariant.  We extend a function $ f $ on $ (t^{-1} F[t^{-1} ]^{2n} \otimes V)_{\bf A}$ to a function 
 on  $S^1 \times  (t^{-1} F[t^{-1} ]^{2n} \otimes V)_{\bf A}$ which we still denote by $f$ by
\[  f ( s , x ) = s f (x ) . \]
Let $ G_r $ and $G_{1 , r } $ denote the isotropy subgroups of $r$ and $ ( 1 , r )$ in $ G( F \la t \ra_+)$, then of course $G_{1 , r }\subset G_r $.
For $ r \in   t^{-1} F[t^{-1} ]^{2n} \otimes V$,  $g \in G_r $, 
\[  \psi ( \frac 12 \la r \alpha_h , r \beta_h \ra ) =  \psi ( \frac 12 \la r , r \beta_h \ra ) = 1 . \]
because $ \psi = 1$ on $F$. Hence 
\[  G_{1 , r } = G_r .\]
Let ${\cal O}$ be a set of representatives
 of $G(F \la t \ra_+)$-orbits in $  tF[t^{-1} ]^{2n}\otimes V$. 
Then

\begin{eqnarray}\label{osum}
    I ( f)  &=&    \int_{ G( F \la t \ra_+ )  \backslash G( {\bf A} \la  t \ra_+ ) } \sum_{r \in t^{-1} F[ t^{-1} ]^{2n} }   f ( (1, r)  h  )dh   \nonumber \\
      &=& \sum_{ r \in {\cal O}}  \int_{ G( F \la t\ra_+ )_r   \backslash G( {\bf A} \la t \ra_+ ) } f ( (1, r) h  )dh .
 \end{eqnarray}

\noindent The classification of $ G( F \la t \ra_+ )$-orbits (same as  $ G( F [[ t ]] )$-orbits)
  in $   t^{-1} F[t^{-1} ]^{2n} \otimes V $ is already done in Section 5 \cite{gz}.
 We recall the result:  if 

\[ x = \sum_{i=1}^l u_i \otimes v_i \in  t^{-1} F[t^{-1} ]^{2n} \otimes V =  t^{-1} F[t^{-1} ]^{2n} \otimes_{F[[ t]] } V[[ t]], \]

\noindent it defines an $F[[t]]$-map 

\[   f_x :   V[[ t]] \to  t^{-1} F[t^{-1} ]^{2n}, \, \, \, \, \, \, \, f_x ( v ) = \sum_i (v_i , v ) u_i  ,\] 

\noindent here $ ( , )$ denotes the $F[[t]]$-valued bilinear form on $V[[ t]]$ that extends  the bilinear form $( , )$ on $V$.
The image $ {\rm Im } f_x $ of $f_x $ is  finite dimensional, so   ${\rm Im }f_x \in Gr ( t^{-1} F[ t^{-1} ]^{2n} ) $ .
 And $x$ can be written as 
 $ x = \sum_{i=1}^l  u_i \otimes v_i $ for a quasi-basis $ u_1 , \dots , u_l $ of ${\rm Im } f_x $ ( Lemma 5.3 \cite{gz}).
We have $ T ( x ) = \sum_{i, j } ( v_i , v_j ) u_i \otimes u_ j \in S_t^2 ( {\rm Im } f_x  )$ (see Section 5 \cite{gz}).
We have 

\begin{theorem}\label{class}
The $ G ( F[[t]])$-orbit in $  t^{-1} F[t^{-1} ]^{2n} \otimes V $ are in one-to-one correspondence with the set of pairs 
 $ W\in Gr ( t^{-1} F[t^{-1}]^{2n} )$, $i \in S_t^2 ( W)$ such that 
  $U(i)_F$ is not empty. The correspondence is that, the orbit of $x$ corresponds to the pair (${\rm Im } f_x $, $T(x)$).  
\end{theorem}

\noindent  This Theorem follows from Theorem 5.8 \cite{gz} directly.  We put 

\[  I_W ( f ) =  \sum_{ r \in {\cal O}: {\rm Im} f_r = W  }  \int_{ G( F [[ t]] )_r   \backslash G( {\bf A} [[ t]] ) } f ( (1, r)  h  )dh  \]

\noindent Therefore 
\begin{equation}\label{7de} I ( f ) = \sum_{ W \in Gr ( t^{-1} F [ t^{-1} ]^{2n} ) }  I_W ( f ) .\end{equation}
 
\begin{lemma}  Let $W\in Gr ( t^{-1} F[ t^{-1}]^{2n} ) $ be  
in the image of  $ P : Gr ( F(( t))^{2n} ) \to Gr ( t^{-1}F[t^{-1}]^{2n})$, 
then 
\begin{equation}\label{7co}   E_W ( f ) =  I_W ( f ) .\end{equation}
 \end{lemma}

\noindent {\it Proof.} We consider the snt-module $M= \tilde W$ as defined in Section 6.
   By Theorem 8.1 \cite{gz},   $ {\rm Et}_W  ( f_1 ) = {\rm It}_W ( f_1 )$ for 
 $f_1 \in {\cal S} ( (W\otimes V)_{\bf A} )$.  Notice that when $ f_1 $ is the restriction 
 of $f$ on $(W\otimes V)_{\bf A}$, we have 
 $ E_W ( f ) = {\rm Et}_W ( f_1 )$, and $ I_W ( f ) = {\rm It}_W ( f_1 )$,
  so (\ref{7co}) holds.
\hfill $\Box$

By the above lemma, and using the fact that $E_W ( f ) = 0 $ for  $W\notin {\rm Im} P $,  we see that 
 \begin{equation}\label {7di}
 I ( f ) - E ( f ) =  \sum_{ W\in Gr ( t^{-1} F[ t^{-1}]^{2n} ) : W\notin {\rm Im} P }  I_W ( f) .
 \end{equation}  

In the Section 8, we shall prove 
\begin{lemma}\label{vanish} If $ W \notin {\rm Im } P $, then 
 $  I_W ( f ) = 0 $.
\end{lemma}

From this lemma, we have our main result:

\begin{theorem}\label{main}
 If ${\rm dim } V > 6n+2 $,    for $ T , g , \phi $ as in Theorem \ref{conv2},  
 \[  E ( Tg \phi   ) = I ( Tg \phi   ) . \]
\end{theorem}

\

\section{Proof of Lemma \ref{vanish}}

In this section we prove Lemma \ref{vanish} that is used to prove Theorem \ref{main}.
  We consider the $F[[t]]$-module 
\[  X =  F((t))^{2n} \otimes_{ F(( t)) } V (( t )) /   F[[ t]]^{2n} \otimes_{F[[t]]} V[[ t ]] .\]
 It is clear that $X$ is isomorphic to the $F[[t]]$-module 
  \[ F(( t )) / F[[t]] \otimes F^{2n} \otimes F((t)) .\] 
 Since  $ {\rm Sp}_{2n} ( F[[t]]) \times G ( F[[t]])$ acts on  $F((t))^{2n} \otimes_{ F(( t)) } V (( t ))$ and 
 preserves the subspace $ F[[ t]]^{2n} \otimes_{F[[t]]} V[[ t ]]$, it acts on 
 $X$ as $F[[t]]$-module isomorphisms. 
 We view $  t^{-1} F [ t^{-1} ]^{2n}$ as an $F[[t]]$-module 
 by the identification 
  \begin{equation}\label{8.1}  t^{-1} F [ t^{-1} ]^{2n} =  F(( t))^{2n} /   F [[  t ]]^{2n} . \end{equation}
Since the $ {\rm Sp}_{2n} ( F[[t]] )$-action on  $ F(( t))^{2n}$ leaves $F [[  t ]]^{2n}$ invariant, 
    it acts on  $t^{-1} F [ t^{-1} ]^{2n}$ as $F[[t]]$-module isomorphisms by the identification (\ref{8.1}). 
  It is clear that 
 \[ X =  F((t))^{2n} \otimes V  /   F[[ t]]^{2n} \otimes V =   t^{-1} F [ t^{-1} ]^{2n} \otimes V .  \]
 \noindent  Similarly, we  view $ t^{-1}V[  t^{-1} ]$ as an $F[[t]]$-module by the identification
 \[  t^{-1} V [ t^{-1} ]  = V (( t))  /   V [[  t ]]  , \]
on which $G ( F[[t]]$ acts as $F[[t]]$-module isomorphisms.
  We see that the following $F[[t]]$-modules

\begin{eqnarray}
   &&t^{-1} F [ t^{-1} ]^{2n} \otimes V,  \, \, \, \, \, \, \, \,   t^{-1} F [ t^{-1} ]^{2n} \otimes_{F[[t]]}V[[ t]] \nonumber  \\
    &&  F^{2n} \otimes t^{-1}  V [ t^{-1} ], \, \, \, \, \, \, \, \,  F[[t]]^{2n}  \otimes_{F[[t]]} t^{-1}  V [ t^{-1} ]. \nonumber
\end{eqnarray}
are all naturally isomorphic to $X$. From now on, we identify the above four spaces.  For example,
        $  e_i t^{-k} \otimes  v$ in the first space is equal to 
        $ e_i \otimes v t^{-k}$ in the third space.
We have the obvious action of   $ {\rm Sp}_{2n} ( F[[t]]) \times G ( F[[t]])$ on
 $t^{-1} F [ t^{-1} ]^{2n} \otimes_{F[[t]]}V[[ t]]$ and $F[[t]]^{2n}  \otimes_{F[[t]]} t^{-1}  V [ t^{-1} ]$.

Every $ x \in  X$, viewed  as an element 
  $ \sum_i u_i \otimes v_i$ in $t^{-1}F [ t^{-1} ]^{2n} \otimes_{F[[t]]}V[[ t]] $ defines an $F[[t]]$-linear map

\[   f_x :   V[[ t]] \to  t^{-1} F[t^{-1} ]^{2n}, \, \, \, \, \, \, \, f_x ( v ) = \sum_i (v_i , v )  u_i  ,\] 

\noindent here $ (v_i , v)$ denotes the $F[[t]]$-valued bilinear form on $V[[ t]]$ that  extends the bilinear form $( , )$ on $V$.
The image $ {\rm Im } f_x $ of $f_x $ is  an $ F[[t]]$-submodule of $ t^{-1} F[ t^{-1} ]^{2n}$, 
 i.e., ${\rm Im } f_x  \in Gr ( t^{-1} F[ t^{-1} ]^{2n} ) $. 
 Note that $x$ can be written as 
 $ x = \sum_{i=1}^l  u_i \otimes v_i $ for a quasi-basis $ u_1 , \dots , u_l $ of ${\rm Im } f_x $ (\cite{gz} , Lemma 5.3).
We have $ T ( x ) = \sum_{i, j } ( v_i , v_j ) u_i \otimes u_ j \in S_t^2 ( {\rm Im } f_x  )$ (see \cite{gz}, Section 5).

 On the other hand,  $x$ defines an $F[[t]]$-linear map:
  \[ f'_x :   F[[ t]]^{2n} \to  t^{-1} V[t^{-1} ], \, \, \, \, \, \, \, f'_x ( v ) = \sum_i \la  a_i , v \ra b_i  ,\] 
 where we write 
 \[  x = \sum a_i \otimes b_i \in F[[t]]^{2n} \otimes_{ F[[ t]] }  t^{-1} V [t^{-1} ] ,\] 
 and $ \la a_i , v \ra$ denotes the $F[[t]]$-valued symplectic form on $ F[[t]]^{2n}$ that extends 
 the standard symplectic form on $ F^{2n}$. We have the following analog of Lemma 5.3 in \cite{gz}:

\begin{lemma}\label{lemma8.1}  Let $x \in X$, and let $ b_1 , \dots , b_m $ be a quasi-basis
  of ${\rm Im } f'_x $, then there are elements $a_1 , \dots , a_m \in F[[t]]^{2n}$ such that 
\newline (1) ${\rm Span}_{ F[[t]]} \{ a_1 , \dots , a_m \}$ is a primitive submodule of $F[[t]]^{2n}$ and 
     $a_1 , \dots ,a_m$ is a basis of ${\rm Span}_{ F[[t]]} \{ a_1 , \dots , a_m \}$.
\newline (2)  $ x = a_1 \otimes b_1 + \dots + a_m \otimes b_m $.
\end{lemma}

The proof is similar to that of  Lemma 5.3 in \cite{gz}, with the role of symplectic space $F^{2n}$ and quadratic space $ V$ interchanged.

For the above $x$, we define 
  \[  T' ( x ) = \sum_{i, j } \la  a_i , a_j \ra b_i \otimes b_j \in \wedge_t^2 ( {\rm Im}  f'_x  ) .\]  
It is clear that for $ g \in {\rm Sp}_{2n} ( F[[ t ]] )$ and $ h \in G ( F[[t]])$, we have 
  \begin{equation}\label{8.2} 
 {\rm Im} f_{ x \cdot h } =  {\rm Im} f_{ x  } , \, \, \, 
  {\rm Im} f'_{ x \cdot g } =  {\rm Im} f'_{ x  } , \, \, \, 
   T (  x \cdot h ) =  T ( x  ) , \, \, \, 
   T' ( x \cdot g ) =  T' ( x ) .
 \end{equation}

\begin{lemma}\label{lemma8.2}
  If $x\in  t^{-1} F [ t^{-1} ]^{2n} \otimes V$, then 
 ${\rm Im } f_x \in {\rm Im } P $ iff $T' ( x ) = 0 $.
\end{lemma}

\noindent{\it Proof. }  If ${\rm Im } f_x \in {\rm Im } P $,  by Lemma 6.2, 
 there is $ g \in {\rm Sp}_{2n} ( F[[t]]) $ such that 
\[     {\rm Im } f_x \cdot g = {\rm Span }_{F[[t]]} ( t^{-k_1}e_1 ,  t^{-k_2}e_2, \dots , t^{-k_l} e_l ) ,\]
where $ l \leq n $.
Using Lemma 5.3 \cite{gz}, we can write 
 \[  x {\cdot g } = e_1 t^{-k_1} \otimes w_1 + e_2 t^{-k_2 } \otimes w_2 + \dots + e_l t^{-k_l } \otimes w_l . \] 
 Note that $ e_i t^{-k_i } \otimes w_i$, under the identification 
    $ X =   F[[t]]^{2n}  \otimes_{F[[t]]} t^{-1}  V [ t^{-1} ]$, is in $ F[[ t]]e_i \otimes_{F[[t]] } t^{-1} V[t^{-1}] $. 
 So $  x {\cdot g }$ is in $\sum_{i=1}^l F[[ t]]e_i \otimes_{F[[t]] } t^{-1} V[t^{-1}] $, 
 therefore $ T' (  x \cdot g ) = 0 $.
By (\ref{8.2}),  
\[ T'( x ) = T' ( x \cdot g ) = 0 .\]
 Conversely, if $T' ( x ) = 0 $.  Let 
 \[ x = a_1 \otimes b_1  + \dots + a_m \otimes b_m \]
 be as in Lemma \ref{lemma8.1}, 
 so 
\[ {\rm Im} f'_x  \cong   F[[t]]/ ( t^{k_1} ) b_1 \oplus \dots   \oplus F[[t]]/ ( t^{k_m} ) b_m  ,\] 
 Then 
 \[ T' ( x ) =\sum_{i, j =1}^m  \la a_i , a_j \ra  b_i \otimes b_j  = 0 \]  
 implies that 
  \begin{equation} \label{8.3} \la a_i , a_j \ra = 0 \, \, \, \, {\rm mod } \, t^{{\rm min} (k_i , k_j) } \end{equation} 
for all $i , j$.  Since $a_1 , \dots , a_m$ satisfies (1) in Lemma \ref{lemma8.1}, we see 
  $a_1  , \dots ,a_m$ project to  linearly independent elements  in $ F[[t]]^{2n} /  t F[[t]]^{2n} $. And (\ref{8.3}) implies that 
 the constant coefficients of $a_1 , \dots ,a_m $ span an isotopic subspace of $F^{2n}$. Therefore   
   $ m \leq n $.   Taking the standard symplectic basis $e_1 , \dots , e_{2n}$ of $F[[t]]^{2n}$, we have 
\[   \la  a_i , a_j \ra = \la e_i , e_j \ra   \, \, \,  {\rm mod } \, t^{{\rm min} (k_i , k_j) } , \]
for $ 1\leq i , j \leq m $. We claim that we can alter $a_1 , \dots , a_m$ to another set
 $ \tilde a_1 , \dots , \tilde a_m$ such that 
 \begin{equation}\label{8.4}
  \tilde a_i = a_i  \, \, \, {\rm mod} \, t^{k_i },  {\rm for } \,   1 \leq i \leq m \end{equation}
and 
\begin{equation}\label{8.5}    \la \tilde a_i , \tilde a_j \ra = \la e_i , e_j \ra = 0   \, \, {\rm for } \,   1 \leq i, j  \leq m .
\end{equation}
The proof of this claim is similar to that of Lemma 5.6 in \cite{gz}, 
 with the role of the symplectic space $F^{2n}$ and quadratic space $V$ interchanged.
  The equation (\ref{8.4}) implies that 
 \[ x = \tilde  a_1 \otimes b_1  + \dots + \tilde a_m \otimes b_m .\]
 In turn  (\ref{8.5}) implies that we can extend the $\tilde a_i$ to a symplectic basis $\tilde a_1 , \dots , \tilde a_{2n}$ of $F[[t]]^{2n}$.
 Then the $F[[t]]$-isomorphism $ g: e_i \mapsto \tilde a_i $ ( $i=1 , \dots , 2n$) is in ${\rm Sp}_{2n} ( F[[ t]] )$.
 So $ x = y \cdot g $, where 
$ y =  e_1 \otimes b_1 + \dots + e_m \otimes b_m $.
 So we have $ {\rm Im}  f_x  = {\rm Im } f_{y}  \cdot g $. It is clear that ${\rm Im } f_y  $ 
   is an $F[[t]]$-submodule of $ W\deff {\rm Span}_{ F[[t]] } \{  e_1 t^{-l_1} , \dots , e_m t^{-l_m} \}$ for non-negative integers $l_1 , \dots , l_n $.
 Since $ m \leq n $,  $ W$  is in ${\rm Im } P$. It follows from Lemma \ref{add} below that $ {\rm Im } f_y$ is in  ${\rm Im } P$.
  Therefore 
${\rm Im}  f_x =  {\rm Im } f_y  \cdot g \in {\rm Im } P $. This completes the proof.
\hfill $\Box $

\begin{lemma}\label{add}
  If $W \in Gr ( t^{-1} F [ t^{-1} ]^{2n } ) $ is in ${\rm Im } P$ and $ W_1 \subset W $ is an $F[[ t]]$-submodule. Then 
  $W_1$ is also in   ${\rm Im } P$.
\end{lemma}

\noindent {\it Proof.} By the assumption,  $W = \pi_- (U)$ for some $U \in Gr ( F(( t))^{2n} )$.
   Let $U'  =  \{  x \in U \, | \,  pi_- ( x ) \in W_1 \}$ and 
 \[ U'' =\{ y \in F[[t]]^{2n} \, | \,  \la y , x \ra = 0 \; {\rm for \; all } \; x \in U' \} . \]  
Then $ U_1 \deff U' + U'' \in Gr ( F(( t))^{2n} ) $ and $ \pi_- ( U_1 ) = W_1 $.
\hfill $\Box $

\

 We consider the group
   \[ O_1 ( F [[ t]] ) \deff \{  g \in O ( F[[t]] ) \, \, | \, \, g \equiv 1 \, \, \, {\rm mod } \, t \} . \]
It is an inverse limit of finite dimensional unipotent groups
    $O_1 ( F [[ t]] ) /  O_k (  F [[ t]])$,  $ O_k ( F[[t]])$ being the subgroup of  $O_1 ( F [[ t]] )$
 consisting of $g$ such that $ g \equiv 1 $ mod $t^k$.
 Let $o_1 ( F [[ t]] )$ denote 
 the space of all $F((t))$-linear maps
\[  \alpha :  V ((t))\to V ((t)) \]
 such that 
 \begin{equation}\label{8.9}  V[[ t]]  \alpha \subset t V[[t]]  \, \, \, \, \, {\rm and } \, \, \, \, \, (v_1 g , v _2 ) + ( v_1 , v_1 g ) = 0                     \end{equation} 
for all $v_1 , v_2 \in V[[t]]$, where $ ( , )$ denote the $F[[t]]$-valued bilinear form on $V[[t]]$ that extends
 $( , )$ on $V$.  Since we assume group $ O(F[[t]])$ acts on $V[[t]]$ from the right, we assume $End_{F[[t]]} V[[t]]$ 
 operates on $V[[t]]$ from the right as well.    The Lie algebra  $o_1 ( F [[ t]] )$ is an inverse limit of finite dimensional
 nilpotent Lie algebras $ o_1 ( F [[ t]] ) /  o_k$, where $o_k$ consists of $\alpha $ such that $ \alpha \equiv 0 $ mod
 $t^k$. We have the exponential map:

  \[ exp :    o_1 ( F [[ t]] ) \to     O_1 ( F [[ t]] ), \, \, \, \, \,  exp( \alpha ) = \sum_{k=0}^{\infty} \frac {\alpha^k } { k! } \] 
 which is a bijection. 

 For $ x\in X$ with $ {\rm Im } f_x  \notin  {\rm Im } P  $, we have $ T' ( x  ) \ne 0 $  by Lemma \ref{lemma8.2}. 
 We let $ W = {\rm Im } f'_x$. Then  $W$ is an $F[[t]]$-submodule of $t^{-1}V[t^{-1}]$.  
    Consider $W + V[[t]] \subset V(( t))$.  It is an $F[[t]]$-submodule of $V(( t))$.  Let 
 \[ \frak g_W = \{    \alpha \in o_1 ( F [[ t]] ) \, \, \, | \, \, \, ( W + V[[ t]] ) \alpha  \subset V[[ t]] \} .\]
 It is easy to see that $ \frak g_W$ is a Lie subalgebra of $ o_1 ( F [[ t]] )$. 
 In $V(( t))$, we define an $F$-valued bilinear form by 
  \[   ( u , v)_F = {\rm Res} ( u , v ) ,  \]
where $( u , v)$ is the $F((t))$-valued bilinear form that extends the bilinear form on $V$, and ${\rm Res} ( u , v )$ is the 
 the coefficient of $t^{-1}$ in $(u, v )$. We have, for $u , v \in V((t))$, 
 \begin{equation}\label{8.5b}
  (   t u , v ) = ( u , t v ) . \end{equation} 
It is easy to see that $V[[t]]$ and $t^{-1} V[ t^{-1} ]$ are maximal 
 isotropic subspaces of $( V((t)), ( , )_F)$. Set 
 \[  W^{\bot } = \{    u \in   V[[ t]] \, | \, ( u , w )_F = 0 \, \, {\rm for \, \, all }\, \, w \in W \}. \]  
Because of (\ref{8.5b}),  $ W^{\bot }$ is an $F[[t]]$-submodule of $V[[t]]$. 
 We let
\[  W' = V[[t]]/ W^{\bot }. \]
We have a non-degenerate pairing $ ( , )_F :  W\times W' \to F$ induced from $( ,)_F$ on $ V((t))$,
 which satisfies (\ref{8.5b}). 
 We call an $F[[t]]$-linear map $ \phi :      W \to W' $  skew symmetric if 
  \[  ( w_1,  \phi ( w_2  ) )_F = -( w_2 ,  \phi  ( w_1 ) )_F . \]  
The space of skew-symmetric $F[[t]]$-linear maps is identified with the dual space of $\wedge^2_t (W )$, the subspace
 of all skew-symmetric tensors in $ W \otimes_{F[[t]]} W$.
 The identification is as follows: the paring of $\phi $ with $ \sum_{i=1}^l w_i \otimes \bar w_i\in \wedge^2_t (W )$ is 
  \[ ( \phi , \sum_{i=1}^l w_i \otimes \bar w_i ) =  \sum_{i=1}^l ( w_i , \phi (  \bar w_i ) ) _F . \]
For each $\alpha \in \frak g_W $, we have 
\begin{equation}\label{8.6}
  V[[t]] \alpha \subset W^{\bot }.
 \end{equation} 
The proof is as follows. For $v\in V[[t]], w\in W$,
\[   ( w , v \alpha )_F = {\rm Res} ( w , v \alpha ) = {\rm Res}( - ( w \alpha , v ) ) = 0 .\] 
The last equality follows from the fact $ w \alpha \in V[[t]]$ and $ ( w \alpha , v )\in F[[t]]$. 
 Each $ g \in \frak g_W $ induces an $F[[t]]$-map
\[  \bar g :  W = (W + V[[t]])  / V[[ t ]] \to  V[[ t]] / W^{\bot } = W' .\]
We have a commutative diagram 
\begin{equation}\label{comm}
\begin{array}{ccc}
  W + V[[t]]   & \overset{g} {\longrightarrow } & V[[ t]]  \\ 
\downarrow \pi   &    & \downarrow \pi'  \\ 
  W  &  \overset {\bar g} {\longrightarrow } &  W' 
\end{array}
\end{equation}
Where $\pi:  W+ V[[t]] \to  W =  W+ V[[t]]/ V[[t]] $ and $\pi':   V[[t]] \to W' = V[[t]] / W^{\bot} $
  denote the natural projections.  Since $ g $ is skew-symmetric, $\bar g$ is skew-symmetric.
 So we have a map 
 \begin{equation}\label{8.7}
  \frak g_W \to ( \wedge_t^2 ( W ))^*  , \, \, \,  g \mapsto \bar g .\end{equation}

\begin{lemma}\label{lemma8a} The above map (\ref{8.7}) is onto.
\end{lemma} 

Before proving this lemma, we need 
\begin{lemma}\label{lemma8b} (1) The rank of every finite dimensional $F[[t]]$-submodule $W$ of $ t^{-1} V [ t^{-1}]$ is at most  
 $ {\rm dim} V$. (2) If the rank of $W$ is ${\rm dim } V$, then $ t^{-1} V[[t]] \subset  W+ V[[t]] $.
 (3). If the rank of $W$ is less than ${\rm dim } V$, we can find an $F[[t]]$-submodule  $W^c$ in $ t^{-1} V [ t^{-1}]$
 such that $ W\cap W^c = \{ 0 \} $ and $ W +  W^c$ has rank ${\rm dim } V $.
\end{lemma} 

\noindent {\it Proof.}   Recall that the rank of $W$ is equal to $ {\rm dim } W/ t W $ and also equal to the number
 of elements in a quasi-basis of $W$.  It is easy to see that $W$ has  a quasi-basis $w_1 , \dots  , w_r$ 
  of the form  \[ w_i = v_i t^{-k_i } + {\rm higher\, terms} \] such that $v_1 , \dots , v_r \in V$ are linearly
 independent. Therefore $ r \leq {\rm dim } V $. This proves (1).   For (2) and  for the quasi-basis $w_1 , \dots , w_r$ as above,
 since $r = {\rm dim } V$,  $v_1, \dots , v_r $ is basis of $V$. Then 
  $ W$ contains $ t^{k_1 -1 } w_1 = t^{-1} v_1 , \dots  , t^{k_r -1 } w_r = t^{-1} v_r $.  Therefore 
 $W+ V[[t]] \supset t^{-1} V[[ t ]] $.  For (3), again we take a quasi-basis $ w_1 , \dots , w_r $ as above,
 since $ v_1 , \dots , v_r$ is linearly independent, we can find  $m = {\rm dim } V - r $ vectors $ v_{1 }^c , \dots , v_{m }^c $
 so that $v_1 , \dots , v_r , v_{1 }^c , \dots , v_{m }^c $ is a basis of $V$. Then 
\[ W^c = {\rm span}_{ F[[ t ]] } \{   t^{-1} v_{1 }^c , \dots , t^{-1} v_{m }^c \} \]
 satisfies the requirements in (3). \hfill $\Box $

\

\noindent {\it Proof of Lemma \ref{lemma8a}. } We first prove the lemma under the condition that the rank of $W$ is equal to ${\rm dim } V$.
 By Lemma \ref{lemma8b} (2),  $ t^{-1} V[[t]] \subset  W+ V[[t]]$.  If $g\in End _{ F[[t]] } ( V(( t)) $ satisfies 
  $ ( W+ V[[t]] ) g \subset  V[[ t ]] $, then $ ( t^{-1} V[[t]] ) g \subset V[[ t]] $, which implies that
  $ V[[ t ]] g \subset t V[[t]] $. For arbitrary skew-symmetric $ \bar g : W \to W'$, we can find an $F[[t]]$-linear map 
  $ \tilde g :   W+V[[ t]] \to V[[ t]] $ such that the diagram (\ref{comm}) is commutative.
   Since $ W+V[[t]]$ and $V[[t]]$ are free $F[[t]]$-modules, $\tilde g$ extends to an $F((t))$-map
  \[ V(( t)) =  (W+V[[ t]])\otimes_{ F[[t]] } F(( t)) \]
 to  \[ V(( t)) = V[[ t]]\otimes_{ F[[t]] } F(( t))\] which we still denote by $\tilde g $. 
Let $ \tilde g^* : V(( t)) \to V(( t)) $ be the dual of $ \tilde g : V(( t)) \to V(( t)) $, i.e., 
 $\tilde g^* $ satisfies 
 \[   ( v_1 \tilde g^* , v_2 ) =  ( v_1  , v_2\tilde g  ) \]
 for all $ v_1 , v_2 \in V(( t)) $.  We claim $ \tilde g- \tilde g^* \in o_1 ( F[[t]] ) $. 
 It is clear that $  \tilde g- \tilde g^*$ satisfies the 2nd condition of (\ref{8.9}). To prove it also satisfies the 1st condition of
  (\ref{8.9}), we need to prove  
\begin{equation}\label{8.11}   (W+V[[ t]]) \tilde g^* \subset  V[[ t]]. \end{equation}
 We take $ w \in W + V[[ t]] , v \in V[[ t]] $, 
 then \begin{equation}\label{8.12} ( w \tilde g^* , v )_F = {\rm Res} ( w \tilde g^* , v ) = {\rm Res} ( w , v \tilde g ). 
 \end{equation} 
       Since $ v \tilde g \in W^{\bot } $,  (\ref{8.12})$= 0$, and this is true for all $v$, so  $w \tilde g^* \in V[[ t]] $.
 This proves (\ref{8.11}), therefore we have proved $\tilde g- \tilde g^* \in o_1 ( F[[t]] ) $.
    It is clear that the image of $\frac 12 (\tilde g- \tilde g^*) $ is $\bar g $.
  
 For the case that the rank of $W$ is smaller than ${\rm dim } V$, we take $W^c$ as in (3) of Lemma \ref{lemma8b}.
 Then 
 \[   (W + W^c )' \cong   W' \oplus (W^c)' .\]
If we have $ \bar g :  W \to W'$ is in $ \wedge_t^2 ( W ) $, we extend $\bar g $ to $\bar {\bar g} : W + W^c  \cong     (W + W^c )' $
 by setting $ ( w + w^c ) \bar {\bar g } = w \bar g $ ,
 then there is $ g \in \frak g_{ W + W^c } $ such that $ g $ maps to  $\bar {\bar g}$. It is clear that $  \frak g_{ W + W^c } \subset   \frak g_{ W}$.
  Then $g \in    \frak g_{ W}$ maps to $\bar g$. 
\hfill $\Box $

  \

\noindent {\it Proof of Lemma \ref{vanish}.}  It is enough to prove that if $x\in t^{-1} F[ t^{-1}]^{2n}\otimes V$
 such that $ W= {\rm Im} f_x \notin {\rm Im} P$, then 
   \begin{equation}\label{last}   \int_{ G( F \la t \ra_+ )_x   \backslash G( {\bf A} \la t \ra_+ ) } f ( (1, x) h  )dh = 0 .
  \end{equation}
The left hand side of (\ref{last}) equals to 
\begin{equation}\label{last2}
 \int_{ G( F \la t \ra_+ )_x   \backslash G( {\bf A} \la t \ra_+ ) }  \psi (\frac 12 \la x \alpha_h , x \beta_h \ra )   f ( x  \alpha_h  )dh.
\end{equation}
If $ h = h_1 h_2$ with $h_1 \in  G( {\bf A} \la t \ra_+ )_x $, then   
  \[ \alpha_{h} =   \alpha_{h_1}   \alpha_{h_2} ,  \, \, \, \, \, \beta_{h} =   \alpha_{h_1}  \beta_{h_2 } +  \beta_{h_1} \delta_{h_2 } .\]
Using $ x \alpha_{h_1} = x $, we have
 \begin{eqnarray}\label{last3}
 && \la x \alpha_h , x \beta_h \ra   \\
  &&  = \la  x \alpha_{h_1}   \alpha_{h_2} , x ( \alpha_{h_1}  \beta_{h_2 } +  \beta_{h_1} \delta_{h_2 } )\ra  \nonumber \\ 
  && =  \la x \alpha_{h_2} , x \beta_{h_1 } \delta_{h_2 } \ra + \la x \alpha_{h_2} , x \beta_{h_2 } \ra . \nonumber 
\end{eqnarray}
Notice that  \[ 
\la x \alpha_g , y \delta_g \ra = \la x , y \ra   \]
 for arbitrary $g \in \widetilde{Sp}_{2N} ( {\bf A}\la t \ra ) $. We therefore have 
\begin{equation}\label{last4} 
 \la x \alpha_h , x \beta_h \ra =  \la x , x \beta_{h_1 } \ra + \la x \alpha_{h_2} , x \beta_{h_2 } \ra .
\end{equation}
  Using (\ref{last4}), the integral (\ref{last2}) equals 
\begin{equation}\label{last6} 
  \int_{    G( F \la t \ra_+ )_x   \backslash G( {\bf A} \la t \ra_+ )_x }  \psi ( \frac 12    \la x , x \beta_{h_1 } \ra ) dh_1 
\end{equation} 
times 
  \[     \int_{    G( {\bf A}  \la t \ra_+ )_x   \backslash G( {\bf A} \la t \ra_+ ) }  \psi ( \frac 12    \la x \alpha_{h_2} , x \beta_{h_2 } \ra )
        f ( x  \alpha_{h_2}  )dh_2 .\]
It is sufficient to prove (\ref{last6}) is $0$. Notice that 
  the map 
\[ \chi:   G( {\bf A} \la t \ra_+ )_x \to S^1 , \, \, \, \, \, h \mapsto  \psi ( \frac 12    \la x , x \beta_{h } \ra ) \]
  is a character of $  G( {\bf A} \la t \ra_+ )_x$ which is trivial on $ G( F [[ t]] )_x $. Therefore to prove  (\ref{last6})=$0$,
 we only need to prove $\chi $ is non-trivial on $ G( {\bf A} \la t \ra_+ )_x $. It is enough to prove that for each place $v$, 
  the character 
 \[ \chi_v :  G( F_v [[ t]] )_x \to S^1 , \, \, \, \, \, h \mapsto  \psi_v ( \frac 12    \la x , x \beta_{h } \ra ) \]
 is non-trivial.   It is sufficient to prove the group homomorphism 
   \[ \phi:   G( F_v [[ t]] )_x \to A( F_v ):     h \mapsto      \la x , x \beta_{h } \ra \]
 is onto, where $A$ denotes the one-dimensional additive group.  The morphism $\phi$ factors through a finite-dimensional quotient 
  of   $G( F_v [[ t]] )_x$.  It is sufficient to prove  the Lie algebra map $ d \phi $ is not trivial.  
 Let $ W = {\rm Im} f'_x $ and $ \frak g_W ( F_v) = \frak g_W \otimes F_v $.  By Lemma \ref{lemma8a}, we have a surjective map
   \begin{equation}\label{8sur} 
  \frak g_W (F_v ) \to \wedge^2_t (  W_v ).  
  \end{equation}
 And by Lemma \ref{lemma8.2}, we have $T' ( x) \ne 0 $ because of our assumption that ${\rm Im } f_x \notin {\rm Im } P $. 
  We may take $ b \in \frak g_W (F_v )$ such that 
   \[  ( T' ( x) , \bar b )      )\ne 0  ,\]
where $\bar b : W \to W'$ is the map induced from $b:   W+ V[[ t]] \to V[[t]]$
 as in (8.9).
 For every $\epsilon \in F_v $, 
\[  g( \epsilon ) = \sum_{k=0}^{\infty}   \frac { \epsilon^k  b^k } {k !} \in G( F_v [[ t]] ) .\]    
Since  \[           (W + V[[t]])_v b \subset  V[[t]]_v ,\] 
  \[   x b \in V[[t]]_v \] 
so $ x  \alpha_{g(\epsilon)} = x $, i.e.,  $g(\epsilon)\in G( F_v [[t]])_x $ for all $\epsilon \in F_v $, . It is easy to see that 
\[    \la x , x \beta_{ g(\epsilon) } \ra \]
is a polynomial of $\epsilon$ whose $\epsilon$-term is 
\[  ( T' ( x) , \bar b ) \ne 0 . \] 
 This proves that $d \phi $ is not $0$ and therefore Lemma 7.4.  \hfill $\Box$.


\begin{thebibliography}{aa}




\bibitem{Garland1}
H. Garland,  Certain Eisenstein series on loop groups: convergence and the constant term,  
 {\em   Algebraic Groups and Arithmetic}  Tata Inst. Fund. Res. Mumbai (2004), 275-319.

\bibitem{Garland2}
H. Garland, Absolute convergence of Eisenstein series on loop groups,
 {\em Duke Math J.} Vol. 135, No. 2 (2006), 203-260.

\bibitem{gz}H. Garland,  Y.Zhu,  On the Siegel-Weil formula for loop groups (I), preprint, arXiv:0812.3236, 2008. Submitted. 

\bibitem{Lag} R. Langlands,  Euler products, Yale University mimeographed notes, 1967


\bibitem{Ma}
H. Matsumoto,  Sur Les Sous-groupes Arithmetiques des groupes semi-simples depolyes,   
 {\em Ann. Sci. Ecole Norm. Sup.}  (4) 2, (1969), 1-62.

\bibitem{S}
R. Steinberg, Lectures on Chevalley groups. Notes prepared by John Faulkner and Robert Wilson. {\em} Yale University, New Haven, Conn., 1968.



\bibitem{W1}
A. Weil, Sur certaines groupes d'oprateurs unitaries, {\em Acta Math. 111} (1964), 143-211.

\bibitem{W2}
A. Weil, Sur la formule de Siegel dans la theorie des groupes classiques,  {\em Acta Math. 113} (1965), 1-88.


\bibitem{Z}
Y. Zhu,  Theta functions and Weil representations of loop symplectic groups, {\em  Duke Math. J}, vol 143 (2008), no. 1, 17-39.


\end{thebibliography}
\end{document}